\documentclass{article}

\newcommand{\R}{I\!\! R}
\newcommand{\C}{I\!\! C}
\newcommand{\Z}{I\!\! Z}

\newcommand{\Q}{I\!\! Q}

\begin{document}

 Tsemo Aristide

C.P. 79067, Gatineau Canada

 J8Y 6V2

tsemo58@yahoo.ca

\bigskip
\bigskip

\centerline{\bf Gerbes, $2$-gerbes and symplectic fibrations.}

\bigskip
\bigskip

\centerline{\bf 1. Introduction.}

\bigskip

A {\bf symplectic fibration} $P\rightarrow N$ is a differentiable
fibration whose typical fiber is the closed connected symplectic
manifold $(F,\omega)$, and such that there exists a trivialization
$(U_i,g_{ij})$, such that $g_{ij}(u)$ is a symplectic automorphism
of the fiber over $u$, endowed with a symplectic structure
$\omega_u$, symplectomorphic to $(F,\omega)$. We suppose that the
cohomology class $[\omega_u]$ of $\omega_u$ is fixed. The theory
of symplectic bundles have been studied by different authors (see
[8], [9], [13], [14]). One purpose of the paper [14] is to
determine whether the structural group of the symplectic bundle
can be reduced to the Hamiltonian group of $(F,\omega)$, that is
whether there exists a symplectic bundle $P'\rightarrow N$
isomorphic to $P$, whose coordinate changes  $g'_{ij}(u)$  are
Hamiltonian automorphisms of the fiber above $u$; such a reduction
will be called a {\bf Hamiltonian structure}, or a
$Ham$-reduction. In [14], it is shown that  the existence of such
Hamiltonian reductions on a finite cover of $N$ is equivalent to
the  following two conditions:

(i) There exists a closed $2$-form $\Omega$ defined on $P$ whose
cohomology class $[\Omega]$   extends $[\omega]$. This means that
the restriction to the fiber above $u$ of the cohomology class
$[\Omega]$, is the cohomology class $[\omega]$. Following McDuff,
we will call the form $\Omega$ a {\bf closed connection form}.

(ii) Let  $Symp(F,\omega)_0$ be  the connected component of the
group of symplectomorphisms $Symp(F,\omega)$, of $(F,\omega)$. The
symplectic bundle is isomorphic to a symplectic bundle whose
coordinate changes take their values in $Symp(F,\omega)_0$.

In [14] it was necessary to impose condition (ii) because the
Hamiltonian subgroup is connected.
  In [14],  McDuff has defined a disconnected subgroup $Ham^s$
of the group  $Symp(F,\omega)$, and has shown that the existence
of a $Ham^s$-reduction of a symplectic bundle is equivalent to the
existence of a closed connection form.

\bigskip

One purpose of this paper is to study the problem  of the
existence of  Hamiltonian and $Ham^s$-reductions of a symplectic
bundle using gerbes, and $2$-gerbes. The theory of gerbes has been
defined by Giraud [6] with the purpose of giving geometric
interpretations of cohomology classes. These classes represent the
obstruction to globally extending locally defined   bundles, as it
is the case for Hamiltonian bundles. Lawrence Breen [2] has also
defined a theory of $2$-gerbes. A $2$-gerbe represents
geometrically the obstruction for a $2$-geometric type structure
to be defined globally. This theory will be also involved here.
For $n\geq 2$, such a geometric obstruction theory has been
defined by Tsemo [20].

Let $\omega$ be a $2$-closed form defined on the manifold $F$, and
$T^1$ the circle. It has been shown by Kostant and Weil, that  the
cohomology class $[\omega]$ of $\omega$ is integral, if and only
if $[\omega]$ is the Chern class of a $T^1$-bundle. When the class
is not integral, we define a flat gerbe $C'(\omega)$ bounded by
the sheaf of locally constant ${\R}$-functions defined on $F$
which represents the obstruction of $[\omega]$ to be zero. We can
construct from this gerbe, another gerbe $C(\omega)$ bounded by
the sheaf of locally constant $T^1$-functions defined on $F$,
which represents the obstruction of $[\omega]$ to be integral (see
2.4). These gerbes are used to study the extension of $[\omega]$.
We have:

\bigskip

{\bf Theorem 2.5.4.}

{\it Let $p:P\rightarrow N$ be a symplectic bundle whose typical
fiber is $(F,\omega)$, there exists a gerbe $C^1_F(\omega)$ whose
classifying cocycle $c_F^1(\omega)$ represents the obstruction of
the symplectic bundle $p$ to have a Hamiltonian reduction.}

\bigskip

To show an analogous theorem for $Ham^s$-reductions, one has to
show first as in  [13], that the automorphisms group of a
$Ham^s$-reduction of a symplectic bundle is independent of the
chosen $Ham^s$-reduction, in order to define the band of the
classifying gerbe. We prove also the following result:

\bigskip

{\bf Theorem 8.2, 8.2.2.}

{\it There exists a $2$-gerbe $C^2_F(\omega)$ whose classifying
cocycle $c^2_F(\omega)$ represents the obstruction of the class
$[\omega]$ to be extended to $P$. The class $[c_F^2(\omega)]$ can
be deduced from $[c^1_F(\omega)]$ as follows: Let $L_1$ and $L_0$
be the respective bands of $C_F^1(\omega)$ and $C_F^2(\omega)$.
There exists an exact sequence of sheaves $1\rightarrow
L_0\rightarrow L'_1\rightarrow L_1\rightarrow 1$, such that the
class $[c^2_F(\omega)]$, is the image of the class
$[c^1_F(\omega)]$ by the connecting morphism
$H^2(N,L_1)\rightarrow H^3(N,L_0)$ of the last exact sequence.
This shows that the existence of a Hamiltonian reduction implies
that the form $\omega$ can be extended to $P$.}

\bigskip

 In $[14]$, McDuff defines a discrete subgroup $H^1(F,P_{\omega})$
  of $H^1(F,{\R})$ and a class in
$H^2(N,H^1(F,P_{\omega}))$  which is the obstruction to have a
$Ham^s$-reduction, that is to obtain a closed connection form. We
show that this last class and $[c^2_F(\omega)]$ are the image of
the Chern class of a $H^1(F,{\R})/H^1(F,P_{\omega})$-principal
bundle by connecting morphisms  related to exact sequences of
sheaves see 8.3.

The holonomy of a connective structure defined on a gerbe,  is the
analogous of the holonomy of a connection. It is used to represent
the action in string theory. We relate the holonomy of the gerbe
$C(\omega)$ to the flux  see 4.

\bigskip

We generalize the methods applied here to solve other geometric
problems, as for example to find a $H$-reduction of a $G$-bundle
such that $G/H$ is a $K(\pi,1)$ space. For this problem, we define
also a gerbe $C_H$ which represents the geometric obstruction to
solve it: More precisely we have:

\bigskip

{\bf Theorem 2.6.3.}

{\it Let $f:P\rightarrow N$ be a $G$-bundle defined on $N$, and
$H$ a subgroup of $G$ such that the right quotient of $G$ by $H$,
$G/H$ is a $K(\pi,1)$ space. Suppose that:

 (i) either the coordinate changes take their values in $Nor(H)$, the
normalizer of $H$ in $G$, this condition is satisfied for
symplectic bundles whose coordinate changes take their values in
the connected component $Symp(F,\omega)_0$ of the group of
symplectic automorphisms $Symp(F,\omega)$.  We consider $G$ to be
$Symp(F,\omega)_0$, and  $H$ to be $Ham(F,\omega)$ the group of
Hamiltonian diffeomorphisms.

 (ii) or $H$ intersects every connected
component of $G$, and there exists a commutative group $L$, a
continue and surjective cocycle $F:G\rightarrow L$, for a
representation $\rho:G\rightarrow L$, such that $\rho(G_0)$ the
image of the connected component $G_0$ of $G$ is the identity of
$L$, and the kernel of $F$ is $H$. Here $L$ is a quotient of a
vector space by a discrete subgroup. This condition is satisfied
if $H$ is the subgroup $Ham^s$, and $G$ is $Symp(F,\omega)$

 Then there exists a gerbe
$C_H$, whose classifying cocycle represents the obstruction to
reduce $G$ to $H$.}

\bigskip

When the gerbe $C_H$ is defined by a cocycle $F$, the classifying
cocycle of this gerbe is the Chern class of a $G/H=L$-bundle.

 Analogues of the gerbe which appears in the last theorem can be constructed
 in more abstract situations: we generalize this construction to the case of
 topoi (elementary topoi). This will perhaps suggest applications
 to algebraic geometry and arithmetic.

\bigskip

 The fact that  the pull-back of a
$Symp(F,\omega)_0$-bundle  endowed with a closed connection form
to a finite cover of the base space has Hamiltonian reductions,
  suggests that the natural category for the study of
Hamiltonian reductions is the etale topos of the base (see 3).

The last part of the paper is devoted to geometric quantization.
We give an extension of the Kostant-Souriau quantization whenever
the class $[\omega]$ is not supposed to be integral, using the
gerbe $C(\omega)$. In particular we obtain the following:

\bigskip

{\bf Theorem.}

{\it Let $(M,\omega)$ be a symplectic manifold, and
$(C^{\infty}(M),\{,\})$  the poisson algebra of $(M,\omega)$.
There exists a preHilbert space $H$, and a representation
$(C^{\infty}(M),\{,\})\rightarrow (Aut(H),[.])$ Where
$(Aut(H),[,])$ is the algebra of operators of  $H$ endowed with
the commutator bracket.}

\bigskip

\bigskip

\medskip

\centerline{\bf Contents.}

\bigskip

1. Introduction.

2. Gerbes theory.

2.1. The classifying cocycle of a gerbe.

2.2. Notations.

2.3. Connective structures on gerbes.

2.4. The gerbe associated to a closed $2$-class.

2.5. Symplectic fibrations and gerbes.

2.6. The McDuff construction of $Ham^s$ and closed extension
forms.

2.7. The universal obstruction of McDuff.

2.8. Generalizations to topoi.

3. The group $Ham^s$ and the etale topos of a manifold.

4. Flux and holonomy of gerbes, the classifying cocycle.

5. A geometric interpretation of a section $H_1(F,{\R})\rightarrow
SH_1(F,{\R})$.

6. Existence of symplectic bundles and gerbes.

7. $2-$gerbes and $2-$gerbed towers.

8. The general case.

9. Quantization of the symplectic gerbe.

\bigskip
\bigskip

 {\bf Acknowledgement.}

The author had the idea to write this paper after a talk of Dusa
McDuff on symplectic fibrations given at Montreal on November
2004. The author want to thank Dusa McDuff who read, correct his
manuscript and help him to improve the presentation of the paper,
and for helpful discussions and comments.

\bigskip
\bigskip

\centerline{\bf  2. Gerbes theory.}

\bigskip

The notion of gerbe has been defined by Giraud [6] to give a
geometric interpretation of $2$-Cech cohomology classes, and to
find obstructions to solve gluing problems. The basic example of a
gerbe is defined as follows: consider a $G$-principal bundle
defined on the manifold $N$, and $1\rightarrow H\rightarrow
G'\rightarrow G\rightarrow 1$ a central extension. The geometric
obstruction to the existence of a $G'$-principal bundle over $N$,
whose quotient by $H$ is the original $G$-bundle is defined by the
classifying cocycle of a gerbe. Gerbe theory also has  a lot of
applications in algebraic geometry. In theoretical physics, a
notion of holonomy of gerbe allows us to represent geometrically
the action in string theory. In this part, we summarize the
results of  gerbe theory used here. We prefer the point of view of
sheaf of categories rather to the one of descent.

\bigskip

{\bf Definition 2.1.}

 Let $N$ be a category. A {\bf sieve $T$} is a subclass of
objects  of $N$, such that if $u$ is an element of $T$, and
 $v\rightarrow u$ an arrow of $N$, then $v$ is an element of $T$.

  Recall that the category $N_u$ is the category whose
 objects are objects $v$ of $N$ such that there exists an arrow
 $h_v:v\rightarrow u$, a morphism between two objects $v$ and $v'$
 of $N_u$ is an arrow $h:v\rightarrow v'$ such that $h_{v'}\circ
 h=h_{v}$.

 A {\bf topology on the category} $N$ is defined as follows: for each
  object $u$ of $N$, there is a family of sieves $J(u)$  of $N_u$, such that:

 (i) If $h:v\rightarrow u$ is an arrow, and $T$ an element of $J(u)$,
 then $T^h=\{v'\in Ob(N): v'\in T, \exists h':v'\rightarrow v\}$ is an
 element of $J(v)$.

 (ii) Suppose that $T$ is a sieve of the sub-category $N_u$ above $u$, if for each map
  $h:v\rightarrow u$,  $T^h$ is an element of $J(v)$, then
$T$ is an element of $J(u)$.

\bigskip

For example, one can define a topology $J$ on the category
$Top(N)$, whose objects are  open sets of a topological manifold
$N$, and morphisms are canonical inclusions as follows: For each
open set $U$ of $N$, an element of $J(U)$ is a sieve of the
category above $U$, which contains a family of open subsets of
$(U_i)_{i\in I}$ of $U$ whose union is $U$.

We will suppose in the sequel that our category is a topos, reader
unfamiliar to this notion make the stronger assumptions that the
category is stable by finite sums, and products, there exists
final and initial objects, and the limits exist and are universal.

We will also suppose that the topology is generated by a covering
family $(u_i)_{i\in I}$, where $u_i$ is an object of $N$. This
means that: for each object $u$, there exists a subset $I_u$
contained in $I$, such that for each $i\in I_u$, there exists a
map $u_i\rightarrow u$ of $N$, the subcategory ${u_{(u_i)_{i\in
I_u}}}$ whose objects are objects $v$ of $N$ such that there
exists a map $v\rightarrow u_i$ , $i\in I_u$ is an element of
$J(u)$. An generating family $(U_i)_{i\in I}$ of a topological
space $N$ generates the topology of the category $Top(N)$.

\bigskip

{\bf Definition 2.2.}

 Let $(N,J)$ be a category $N$ endowed with a topology $J$.
  A {\bf sheaf of categories} defined on $(N,J)$ is a
 correspondence $C$:

 $$
U\longrightarrow C(U)
$$

where $C(U)$ is a category, and $U$ an object of $N$, which
verifies the following properties:

(i) For each  map $U\rightarrow V$, there exists a restriction
map $ r_{U,V}: C(V)\longrightarrow C(U)$

such that

$$
r_{U_1,U_2}\circ r_{U_2,U_3}=r_{U_1,U_3}
$$

In fact, while the previous equality is verified in many examples,
only an isomorphism between $r_{U_1,U_2}\circ r_{U_2,U_3}$ and
$r_{U_1,U_3}$ is needed. The last relation defined the notion of
{\bf presheaf of categories.}

 The following properties needed to
be verified to complete the notion of sheaf of categories.

\bigskip

(ii) {\bf Gluing properties for objects.}

\bigskip

Let $(U_i)_{i\in I}$ be a covering  family of the object $U$ of
$N$, and $e_i$ an object of $C(U_i)$. We denote abusively by $N$
the final object of $N$. Suppose there are morphisms

$$
g_{ij}:r_{U_i\times_U U_j,U_j}(e_j)\longrightarrow r_{U_i\times_U
U_j,U_i}(e_i)
$$

 such that on $U_{i_1}\times_NU_{i_2}\times_NU_{i_3}$, the restrictions of the morphisms
$g_{i_1i_2}g_{i_2i_3}$, and $g_{i_1i_3}$ between the respective
restrictions of $e_{i_3}$ and $e_{i_1}$ to
$U_{i_1}\times_NU_{i_2}\times_NU_{i_3}$ are equal. Then there
exists an object $e_U$ of $U$ such that $r_{U_i,U}(e_U)=e_i$.

\bigskip

(iii) {\bf Gluing conditions for maps.}

\bigskip

For each objects $e$, $e'$ of $C(U)$, the correspondence defined
on the category above $U$ by

$$
V\longrightarrow Hom(r_{U,V}(e),r_{U,V}(e'))
$$

is a sheaf of sets.

\bigskip

A correspondence $C$ which satisfies properties $(i)$, $(ii)$ and
$(ii)$ is a { sheaf of categories}. A {\bf gerbe} is a sheaf of
categories which satisfies the following conditions:

(iv) There exists a covering family  $(U_i)_{i\in I}$ of $N$ such
that $C(U_i)$ is not empty for each $i$,

(v) {\bf Local connectivity.}

 For each object $U$ of $N$, there exists a covering family
  $(U_i)_{i\in I}$ of $U$ such that, for each pair of elements $e$ and $e'$ of
$C(U)$, $r_{U_i,U}(e)$ and $r_{U_i,U}(e')$ are isomorphic.

(vi) There exists a sheaf $L$ on $N$ such that for each object
$e_U$ of $C(U)$, $Hom(e_U,e_U)=L(U)$, and this identification
commutes with restrictions an arrows. The sheaf  $L$  is called
{\bf the band} of the gerbe $C$, or we say that the gerbe $C$ is
{\bf bounded} by $L$ $\bullet$

\bigskip

{\bf  The classifying cocycle of a gerbe.}

\bigskip

Let $(U_i)_{i\in I}$ be a covering family of $N$ such that for
each $i$, $C(U_i)$ is not empty, and $e_i$ is an object of
$C(U_i)$. Choose  maps $ g_{ij}:r_{U_i\times_N
U_j,U_j}(e_j)\longrightarrow r_{U_i\times_N U_j,U_i}(e_i)$ for all
$i,j$. Denote by ${g_{i_1i_2}}^{i_3}$ the restriction of
$g_{i_1i_2}$ between the restrictions of $e_{i_2}$ and $e_{i_1}$
to $U_{i_1}\times_NU_{i_2}\times_NU_{i_3}$. Then the map

$$
c_{i_1i_2i_3}={g_{i_1i_2}}^{i_3}{g_{i_2i_3}}^{i_1}{g_{i_3i_1}}^{i_2}
$$

 is an automorphism
of $r_{U_{i_1}\times_NU_{i_2}\times_NU_{i_3},U_{i_1}}(e_1)$ which
may be thought of as an element of
$L(U_{i_1}\times_NU_{i_2}\times_NU_{i_3})$. The assignment
$U_{i_1}\times_N U_{i_2}\times_N U_{i_3}\rightarrow c_{i_1i_2i_3}$
is called the classifying cocycle of the gerbe. If the band is
commutative, it is a Cech-cocycle in the classical sense. It has
been shown by Giraud [6] that the isomorphism classes of gerbes
bounded by the sheaf $L$ is one to one with the Cech cohomology
group $H^2(N,L)$, when $L$ is commutative. If the band is not
commutative $H^2(N,L)$ is defined to be set the of equivalence
classes of gerbes bounded by $L$. The trivial gerbe is a gerbe
such that $C(N)$ is not empty. The elements of $C(N)$ are called
global sections. They are one to one with $H^1(N,L)$.

\bigskip

{\bf 2.2 Notations.}

\bigskip

Let $U_{i_1},...,U_{i_p}$ be objects  of a topos  $N$, and $C$ a
presheaf defined on $N$. We will denote by $U_{i_1..i_p}$ the
fiber product of $U_{i_1}$,...,$U_{i_p}$ on the final object. If
$e_{i_1}$ is an object of $C(U_{i_1})$, ${e_{i_1}}^{i_2...i_p}$
will be the restriction of $e_{i_1}$ to $U_{i_1...i_p}$. For a map
$h:e\rightarrow e'$ between two objects of $C(U_{i_1..i_p})$, we
denote by $h^{i_{p+1..i_n}}$ the restriction of $h$ to a morphism
between $e^{i_{p+1}...i_n}\rightarrow {e'}^{i_{p+1}...i_n}$.

\bigskip

Now we provide details on the classic  example of sheaf of
categories given at the beginning. Consider an extension:

$$
1\longrightarrow H\longrightarrow G'\longrightarrow
G\longrightarrow 1
$$

such that $H$ is a central group in $G'$, and the map
$G'\rightarrow G$ has local sections. Supposed defined a
$G$-principal bundle $p_G$, over $N$. The obstruction of the
existence of a $G'$-principal bundle over $N$ whose quotient by
$H$ is $p_G$, is the cohomology class of the classifying cocycle
of the following gerbe $C_H$ defined on the categories of open
subsets of $N$ as follows: for each subset $U$ of $N$, we define
$C_H(U)$ to be the category whose objects are principal
$G'$-bundles over $U$ whose quotient by $H$ is the restriction of
$p_G$ to $U$. To explicit the classifying cocycle $c_H$ of this
gerbe, consider an open covering $(U_i)_{i\in I}$ of $N$, who
trivializes the bundle $p_G$. We denote by $g_{ij}:U_i\cap
U_j\rightarrow G$ the transition functions. Since the projection
$G'\rightarrow G$ has local sections, we can suppose that we can
lift each map $g_{ij}$ to a map $\hat g_{ij}:U_i\cap
U_j\rightarrow G'$. The classifying cocycle of $C_H$ is defined
by:

$$
c_{i_1i_2i_3}={\hat g_{i_1i_2}}^{i_3}{\hat g_{i_2i_3}}^{i_1}{\hat
g_{i_3i_1}}^{i_2}
$$

This situation applies to the case where $H$ is ${\Z}/2$, $G'$ the
spin group, and $G$ the ortoghonal group $O(n)$. The $O(n)$-bundle
is the orthogonal reduction of the bundle of linear frames of the
$n$-dimensional manifold $N$, defined by a riemannian metric. The
gerbe represents the geometric obstruction of the existence of a
spin structure on $N$. The cocycle in this case is the second
Stiefel-Whitney class.

\bigskip

{\bf 2.3 Connective structures on gerbes.}

\bigskip

The notion of a connective structure on a gerbe has been defined
by Deligne see [3]. It is the analogous to the notion of a
connection on a principal bundle.

\bigskip

{\bf Definition 2.3.1.}

  Consider a gerbe $C$ defined on a manifold whose band is $L$.
  {\bf A connective structure} on $C$, is a correspondence which
associates to each object $e_U$ of $C(U)$ a torsor $Co(e_U)$,
called the torsor of connections, that is an affine space whose
underlying vector space is a subset of the set of $1$-forms
defined on $U$.  The following properties are supposed to be
satisfied by this assignment:

(i)- The correspondence $e_U\rightarrow Co(e_U)$ is functorial
with respect to restrictions  to smaller subsets.

(ii)- For every isomorphism $h:e_U\rightarrow e'_U$ between
objects of $C(U)$, there exists an isomorphism of torsors
$h^*:Co(e_U)\rightarrow Co(e'_U)$ compatible with the composition
of morphisms of $C(U)$, and the restrictions to smaller subsets.

\bigskip

Suppose now that  the band of the gerbe is a $T^1$-sheaf, where
$T^1$ is the circle. Then for each morphism $g$ of the object
$e_U$ of $C$, and $\nabla_{e_U}$ a connection of $Co(e_U)$,

$$
g^*\nabla_{e_U}=\nabla_{e_U}+g^{-1}dg
$$

{\bf A curving} of a connective structure $Co$ is an assignment to
each object $e_U$, and  each element $\nabla$ of $Co(e_U)$, a
$2$-form $D(e_U,\nabla)$ defined on $U$ such that

for each morphism $h:e'_U\rightarrow e_U$,
$D(e_U,\nabla)=D(e'_U,h^*\nabla)$.

If $\alpha$ is a $1$-form on $U$ such that $\nabla+\alpha$ is an
element of $Co(e_U)$, then

$$
D(e_U,\nabla+\alpha)=D(e_U,\nabla)+d\alpha
$$

The assignment $e_U\rightarrow D(e_U,\nabla)$ is compatible with
the restrictions to smaller subsets.

{\bf The curvature of the curving} is the form whose restriction
to each open subset such that $C(U)$ is not empty is
$dD(e_U,\nabla)$, where $e_U$ is an object of $C(U)$, and $\nabla$
an element of $Co(e_U)$ $\bullet$

\bigskip

{\bf 2.4 The gerbe associated to a closed $2$-form.}

\bigskip

Let $(N,\omega)$ be a manifold $N$, endowed with a closed $2$-form
$\omega$; $(N,\omega)$ is often called a {\bf Dirac manifold}.
There exists a Cech-Weil isomorphism between the De Rham
cohomology groups of $N$, and the Cech-cohomology groups of the
sheaf of locally constant ${\R}$-functions defined on $N$. Thus
using the theorem of Giraud [6], we deduce that the cohomology
class $[\omega]$ of $\omega$   classifies   a gerbe $C'(\omega)$
defined on $N$ and bounded by the sheaf of locally constant
${\R}$-functions.

In this part, we present the construction of the classifying
cocycle of this gerbe. This is in fact the classic explanation of
the Cech-Weil isomorphism. This gerbe is the fundamental gerbe
used to define many of the geometric obstructions involved in this
paper.

\bigskip

Let $N$ be a manifold,  $\omega$ a closed $2$-form defined on $N$,
and $(U_i)_{i\in I}$ a cover of $N$ by contractible open subsets.
Without loss of generality, we can suppose that $U_i\cap U_j$ is
connected. The Poincare Lemma implies the existence of a family of
$1$-forms $(\alpha_i)_{i\in I}$ such that

$$
d(\alpha_i)=\omega_{\mid U_i},
$$

 where $\omega_{\mid U_i}$ is the restriction of $\omega$
 to $U_i$. Let $\alpha^i_j$ and $\alpha^j_i$ be the respective restrictions of
 $\alpha_j$ and $\alpha_i$ to $U_i\cap U_j$. Denote by
 $\alpha_{ij}$, the form
$\alpha^i_j-\alpha^j_i$ on $U_i\cap U_j$. The form $\alpha_{ij}$
is closed. By applying the Poincare lemma to $\alpha_{ij}$, we
obtain a family of real  valued functions $u_{ij}$ defined on
$U_i\cap U_j$ such that

$$
d(u_{ij})=\alpha_{ij}.
$$

On $U_{i_1i_2i_3}$, the differential of
$c_{i_1i_2i_3}=u_{i_2i_3}-u_{i_1i_3}+u_{i_1i_2}$ is zero. This
implies that it is a constant map. The family of functions
$c_{i_1i_2i_3}$ is a $2$-Cech cocycle for the sheaf of locally
constant ${\R}$-functions.

 If $c_{i_1i_2i_3}\in{\Z}$, the functions
$h_{ij}=exp(2i\pi u_{ij})$ defines a line bundle over $N$. This
bundle is the well-known Kostant-Weil construction. In this case,
the cohomology class $[\omega]$ of $\omega$ is an element of
$H^2(N,{\Z})$.

Suppose that $[\omega]$ is not necessarily an element of
$H^2(N,{\Z})$. Using   Giraud's theorem concerning the
classification of gerbes, we can associate to $\omega$ a gerbe
$C'(\omega)$ bounded by the sheaf of locally constant
${\R}$-functions, whose classifying cohomology class is the image
of $[\omega]$ by the De Rham Cech isomorphism. This gerbe
represents the  obstruction of the class $[\omega]$ to be zero.
The objects of $C'(\omega)(U)$ when it is not empty, can be
represented by flat ${\R}$-bundles by using the reconstruction
theorem of Giraud presented in Brylinsky [3]. We denote by
$c'_{\omega}$ the classifying cocycle of $C'(\omega)$.

The following proposition describes a gerbe bounded by $T^1$ which
will play a fundamental role in this paper.

\bigskip

{\bf Definition-Proposition 2.4.1.}

{\it

Let $U$ be an open subset of $N$, and denote by $C(\omega)(U)$ the
category whose objects are circle bundles over $U$, endowed with a
connection whose curvature is $\omega_{\mid U}$ the restriction of
$\omega$ to $U$. We will denote by $(e_U,\nabla_{e_U})$ an object
of $C(\omega)(U)$; $e_U$ represents a  $T^1$-bundle, and
$\nabla_{e_U}$ the connection on $e_U$ whose curvature is the
restriction of $\omega$ to $U$. The set of morphisms between two
objects $(e_U,\nabla_{e_U})$, and $(e'_U,\nabla_{e'_U})$ is the
set of morphisms of  differential bundles over the identity
$f:e_U\rightarrow e'_U$ such that
$f^*(\nabla_{e'_U})=\nabla_{e_U}$. The correspondence
$U\rightarrow C(\omega)(U)$ is a gerbe bounded by the sheaf of
locally constant $T^1$-valued functions. The class of its
classifying cocycle is the obstruction of $[\omega]$ to be
integral.}

\bigskip

{\bf Proof.}

 First, we show that $C(\omega)$ is a sheaf of categories.

\bigskip

Gluing conditions for objects:

\bigskip

 Let $(U_i)_{i\in I}$ be an open
cover of an open set $U$ of $N$,  $(e_i,\nabla_{e_i})$ an object
of $C(\omega)(U_i)$, and $g_{ij}:e^i_j\rightarrow e^j_i$ a
morphism such that on $U_{i_1i_2i_3}$,
$g^{i_3}_{i_1i_2}g^{i_1}_{i_2i_3}=g^{i_2}_{i_1i_3}$. Since the
elements of the family $(e_i)_{i\in I}$ are bundles, we deduce
that there exists a bundle $e$ over $U$ whose restriction to $U_i$
is $e_i$. The bundle $e$ is endowed with a connection whose
curvature is the restriction of $\omega$ to $U$ since the
restriction of this curvature to $U_i$ is the restriction of
$\omega$ to $U_i$.

\bigskip

Gluing conditions for arrows:

\bigskip

Let $e,e'$ be a pair of elements of $C(\omega)(U)$, the
correspondence defined on the category of open subsets of $U$ by
$V\rightarrow Hom(r_{U,V}(e),r_{U,V}(e'))$ is a sheaf of sets,
since it is the sheaf of morphisms between two bundles.

It remains to verified that the sheaf of categories is a gerbe.

Let $(U_i)_{i\in I}$ be an open covering of $N$ by contractible
open subsets. For each pair of objects $(e,\nabla_e)$ and
$(e',\nabla_{e'})$ of $C(\omega)(U)$ we have to show that these
objects  are locally isomorphic.

To show this, consider two objects $(e_i,\nabla_{e_i})$ and
$(e'_i,\nabla_{e'_i})$ of $C(\omega)(U_i)$. The bundle $e_i$ and
$e'_i$ are isomorphic to the trivial bundle $U_i\times T^1$. Let
$d$ be the differential, $\nabla_{e_i}=d+\alpha_i$, and
$\nabla_{e'_i}=d+\alpha'_i$. For each section $u:U_i\rightarrow
i{\R}$ of the Lie algebra bundle associated to this  bundle, and
each automorphism $g$ defined by a differentiable map
$U_i\rightarrow T^1$, we have:

$$
g^*(d+\alpha_i)(u)=g^{-1}(d+\alpha_i)(gu))=(g^{-1}dg+d+\alpha_i)(u)\leqno
(1)
$$

Since the connections $\nabla_{e_i}$ and $\nabla_{e'_i}$ have the
same curvature, there exists a function $v_i$ such that
$\alpha'_i=\alpha_i+dv_i$.  We can suppose (or shrinking $U_i$ if
needed)  that the logarithm is defined on $U_i$, thus
$g^{-1}dg=dlog(g)$.  If we take $g=exp(v_i)$,where
$\alpha'_i=\alpha_i+dv_i$, then $g^*(d+\alpha_i)=d+\alpha'_i$. We
obtain that the respective restrictions $(e_i,\nabla_{e_i})$ and
$(e'_i,\nabla_{e'_i})$ of $(e,\nabla_{e})$ and $(e',\nabla_{e'})$
to $U_i$ are isomorphic.

\bigskip

The automorphism group of the object $(e,\nabla_e)$ of
$C(\omega)(U)$ is the group of gauge transformations which
preserve the connection $\nabla_e$. These gauge transformations
are necessarily constant maps, as is shown by $(1)$.

Now we have to interpret geometrically the vanishing of the
cohomology class $[c_{\omega}]$, of the classifying cocycle
$c_{\omega}$ of $C({\omega})$. The theorem of Giraud [6] implies
that this is equivalent to the existence of a global object of the
gerbe, that is a $T^1$-bundle over $N$ whose curvature is
$\omega$. The Kostant-Weil construction implies that this is
equivalent to the fact that the class $[\omega]$ of $\omega$ is
integral $\bullet$.

\bigskip

Now we  establish the relation between the gerbes $C'(\omega)$ and
$C(\omega)$.

\bigskip

{\bf Proposition 2.4.2.}

{\it

  Consider the exact sequence of sheaves of locally constant
  functions:

$$
1\longrightarrow{\Z}\longrightarrow {\R}\longrightarrow
T^1\longrightarrow 1 \leqno(1)
$$

where the map ${\Z}\rightarrow {\R}$ is the canonical injection,
and ${\R}\rightarrow T^1$ is the exponential map of the Lie group
$T^1$, that is the composition of the multiplication by $2\pi i$
and the usual exponential. We obtain
 the following exact sequence in cohomology:

$$
H^1(N,T^1)\rightarrow H^2(N,{\Z})\rightarrow
H^2(N,{\R})\rightarrow H^2(N,T^1)...
$$

The class $[c_{\omega}]$ is the image of the class
$[c'_{\omega}]$, by the map $H^2(N,{\R})\rightarrow H^2(N,T^1)$.}

\bigskip

{\bf Proof.}

Consider an open covering $(U_i)_{i\in I}$ of $N$, such that for
each $i$, $U_i$ is contractible and $U_{i_1...i_p}$ is connected
(using a theorem of Weil, we can suppose $U_{i_1...i_p}$ to be
connected). Let $c_{i_1i_2i_3}$ be the classifying cocycle of the
gerbe $C'_{\omega}$. The image of $[c'_{\omega}]$ by the map
$H^2(N,{\R})\rightarrow H^2(N,T^1)$ is represented by the cocycle
$exp(2i\pi c_{i_1i_2i_3})$. Recall that to construct the cocycle
$c_{i_1i_2i_3}$ we have considered the restriction $\omega_{\mid
U_i}$ of $\omega$ to $U_i$. There exists a form $\alpha_i$ such
that $d\alpha_i=\omega_{\mid U_i}$. We can define the object
$e_i=(U_i\times T^1,d+\alpha_i)$ of $C(\omega)(U_i)$. Let
$\alpha^i_j$ be the restriction of $\alpha_j$ to $U_{ij}$, then
there exists a function $u_{ij}$ such that
$d(u_{ij})=\alpha^i_j-\alpha^j_i$. The functions $exp(2i\pi
u_{ij})$ defines a morphism between $e_j^i$ and $e_i^j$. The
classifying cocycle of $C'(\omega)$ is
$c_{i_1i_2i_3}=u_{i_2i_3}-u_{i_1i_3}+u_{i_1i_2}$, and the
classifying cocycle of $C(\omega)$ is $exp(2i\pi
u_{i_2i_3})exp(-2i\pi u_{i_1i_3})exp(2i\pi u_{i_2i_3})=exp(2i\pi
c_{i_1i_2i_3})$.

\bigskip

Now, we are going to endow the gerbe $C(\omega)$ with a connective
structure.

\bigskip

{\bf Proposition 2.4.3.}

{\it For each open set $U$ of $N$, and the object $e_U$ of
$C(\omega)(U)$, the set $Co(\omega)(e_U)$ of connections defined
on $e_U$ whose curvature is the restriction of $\omega$ to $U$
defines a connective structure on $C(\omega)$. The restriction
$\omega_{\mid U}$ of $\omega$ to $U$, is the curving of each
object $e_U$ of $C(\omega)(U)$. The curvature of this curving is
zero.}

\bigskip

{\bf Proof.}

Let $\alpha$ and $\alpha'$ be two elements of $Co(e_U)$, and
$(U_i)_{i\in I}$ a contractible open cover of $U$. It is a
well-known fact that there exists a $1$-form $v$ such that
$\alpha'=\alpha+v$. The restriction of $e_U$ to $U_i$ is
diffeomorphic to the trivial $T^1$-bundle. This implies that under
this identification, the respective restrictions $\alpha_i$ and
$\alpha'_i$ of the connections $\alpha$ and $\alpha'$ to $U_i$,
have the form $d+u_i$, and $d+u_i+v_{\mid U_i}$ where $u_i$ is a
$1$-form defined on $U_i$, and $v_{\mid U_i}$ is the restriction
of $v$ to $U_i$. The respective curvatures of $d+u_i$ and
$d+u_i+v_{\mid U_i}$ are the $2$-forms $du_i$ and $d(u_i+v)$.
Since they coincide with the restriction of $\omega$ to $U_i$, we
deduce that $dv=0$, thus $Co(\omega)(e_U)$ is an affine space
whose underlying vector space is the vector space of closed
$1$-forms. We deduce that it is a torsor.

The fact that for each automorphism $g$ of $e_U$,
$g^*\nabla_{e_U}=\nabla_{e_U}+g^{-1}dg$ results from the fact that
$\nabla_{e_U}$ is a connection.

For each map $h:e_U\rightarrow e'_U$, we define the map
$h^*:Co(\omega)(e_U)\rightarrow Co(\omega)(e'_U)$, to be the
pull-back of connections by $h^{-1}$. This implies that $h^*$
behave naturally in respect to restrictions to smaller subsets and
compositions.

\bigskip

Let $\nabla_{e_U}$ be an element of $Co(e_U)$, the curvature of
$\nabla_{e_U}$ is the restriction of $\omega$ to $U$,
$\omega_{\mid U}$. It is also the curvature of
${h^{-1}}^*(\nabla{e_U})$. This can be shown using $(1)$. This
implies that $\omega$ defines a curving for this connective
structure. The fact that the curvature of this connective
structure is zero follows from the fact that the form $\omega$ is
closed $\bullet$

\bigskip

 At the end of this paper, we will present
a quantization of symplectic manifolds using the gerbe
$C(\omega)$. This gerbe thus appears to be fundamental in
symplectic geometry.

\bigskip

{\bf 2.5 Symplectic fibrations and gerbes.}

\bigskip

Let $p:P\rightarrow N$ be a symplectic fibration, whose fiber $F$
is the closed symplectic manifold $(F,\omega)$. We study the
following problem: extend $[\omega]$ to a class $[\Omega]$ defined
on $P$, that is, find a cohomology class $[\Omega]\in H^2(P,{\R})$
such that for every $u\in N$, consider the canonical embedding
$i_u:F\rightarrow F_u\rightarrow P$, $i_u^*([\Omega])=[\omega]$. A
result of Thurston [16] implies that in this situation there
exists a form $\Omega$ such that ${i_u}^*\Omega=\omega_u$ for all
$u\in N$.

To use the theory of gerbes, we must suppose that the class
$[\omega]$ of the symplectic form $\omega$ is integral. Thus it is
the Chern class of a circle bundle $h_F$ over $F$.  In the general
case, we will use the gerbe $C'(\omega)$ to define a $2$-gerbe
which represents the geometric obstruction of the class $[\omega]$
to be lifted to $P$. We have the following proposition:

\bigskip

{\bf Proposition. 2.5.1.}

{\it Suppose that $[\omega]$ is integral, and consider for each
open set $U$ of $N$ the category $C_F(\omega)(p^{-1}(U))$  of
circle bundles over $p^{-1}(U)$ whose  Chern class is
$[\Omega_U]$, an element of $H^2(p^{-1}(U),{\R})$ which extends
$[\omega]$. The correspondence defined on the category of open
subsets of $P$ by $p^{-1}(U)\rightarrow C_F(p^{-1}(U))$, defines a
gerbe on $P$, where $P$ is endowed with the topology structure
generated by $p^{-1}(U)$, where $U$ is an open subset of $N$, and
its differential structure is modelled on ${\R}^n\times F$, where
$n$ is the dimension of $N$. The cohomology class of the
classifying cocycle of this gerbe is the obstruction to extend
$[\omega]$.}

\bigskip

{\bf Proof.}

\bigskip

Gluing conditions for objects.

\bigskip

Recall that for every objects $e_U$ and $e'_U$ of
$C_F(\omega)(U)$, $Hom(e_U,e'_U)$ are morphisms of circle bundles
which project to the identity. Let $(U_i)_{i\in I}$ be an open
covering of the open set $U$ of $N$ by open subsets, $e_i$ an
object of $C_F(\omega)(p^{-1}(U_i))$, and $u_{ij}:e_j^i\rightarrow
e_i^j$ a morphism which verifies
${u_{i_1i_2}}^{i_3}{u_{i_2i_3}}^{i_1}={u_{i_1i_3}}^{i_2}$, then
there exists a bundle $e_U$ on $p^{-1}(U)$ whose restriction to
each $U_i$ is $e_i$. This is deduced from the classical definition
of a $T^1$-bundle $e_U$ over $p^{-1}(U)$. Consider a $2$-closed
form $\Omega$ which represents the Chern class of $e_U$, since the
restriction of $e_U$ to $p^{-1}(U_i)$ is $e_i$, its Chern class
which is the restriction of the class $[\Omega]$ of $\Omega$ to
$p^{-1}(U_i)$ is the Chern class of $e_i$. This implies that
$[\Omega]$ extends to $p^{-1}(U)$ the class $[\omega]$.

\bigskip

Gluing condition for arrows

The correspondence defined on the category of open subsets of $U$
by $V\rightarrow Hom({e_U}_{\mid V},{e'_U}_{\mid V})$ defines a
sheaf on this category, since it is a sheaf of morphisms between
two bundles.

This shows that $C_F(\omega)$ is a sheaf of categories. It remains
to prove that it is a gerbe.

\bigskip

Let $(U_i)_{i\in I}$ be a cover of $N$ by contractible open
subsets, $p^{-1}(U_i)=U_i\times F$ this implies that
$H^*(p^{-1}(U_i))=H^*(F)$, thus there exists a class $[\Omega_i]$
on $p^{-1}(U_i)$ which extends $[\omega]$, and which is integral.
Thus $C_F(\omega)(p^{-1}(U_i))$ is not empty.

We deduce from $(1)$ that the group of automorphisms of the
objects of $C_F(\omega)(p^{-1}(U))$ are sections of the sheaf
circle valued functions defined on $p^{-1}(U)$.

Connectivity.

Let $e_U$ and $e'_U$ be a pair of objects of $C_F(\omega)(U)$.
Denote respectively by $e_i$ and $e'_i$ the respective
restrictions of $e_U$ and $e'_U$ to $p^{-1}(U_i)$,where
$(U_i)_{i\in I}$ is an open cover of $U$ by contractible open
subsets. Since $U_i$ is contractible, the Chern class of the
differentiable bundle $e_i$ and $e'_i$ are mapped to $[\omega]$ by
the isomorphism $H^2(U_i\times F,{\R})\rightarrow H^2(F,{\R})$.
This implies they are isomorphic since they have the same Chern
class.

\bigskip

 If the classifying cocycle of the gerbe $C_F(\omega)$ has a trivial
cohomology class, then by a theorem of Giraud [6], the gerbe
$C_F(\omega)$ has a global section $e$. Let $u$ be an element of
the contractible open subset $U_i$ of $N$. The restriction of $e$
to $p^{-1}(U_i)$ is an element $e_i$ of $C_F(\omega)(U_i)$, by
definition, its restriction to $F_u$ has Chern class $[\omega]$
$\bullet$

\bigskip

{\bf Remark.}

 Denote the classifying cocycle of the  gerbe $C_F(\omega)$ by
 $c_F(\omega)$. Its cohomology class is an element of the sheaf
 cohomology group of differentiable functions
$H^2(P,T_1)$. We can consider the exact sequence of sheaves of
differentiable functions:

$$
1\longrightarrow{\Z}\longrightarrow {\R}\longrightarrow
T^1\longrightarrow 1. \leqno(1)
$$

We deduce an isomorphism between $H^2(P,T^1)$ and $H^3(P,{\Z})$,
since $H^*(P,{\R})$ the cohomology of the sheaf of
${\R}$-differentiable functions is zero, because there exist
partitions of unity. Thus the gerbe $C_F(\omega)$ is classified by
an element of $H^3(P,{\Z})$.

\bigskip

In [3] Brylinski has studied the following  problem: Suppose
defined on $P$ a $2$-form $\Omega$ whose restriction to each fiber
is closed, integral and symplectic, find obstructions to build a
closed $2$-form whose restriction to a fiber $F_u$ above $u$,
coincides with the restriction of $\Omega$ on the fiber $F_u$. If
$H^1(F,{\R})=0$, the obstruction to find such a class is a gerbe
$C_p(\omega)$ defined on $N$.

Recall the construction of $C_p(\omega)$. For every open set $U$,
$C_p(\omega)(U)$ is the category whose objects are ${T^1}$-bundles
over $p^{-1}(U)$, endowed with a connection such that the
restriction of its curvature to a fiber $F_u$ above $u$, coincide
with the restriction of $\Omega$ to $F_u$. A morphism between two
objects $(e_U,\nabla_{e_U})$ and $(e'_U,\nabla_{e'_U})$ is a
morphism $f:e_U\rightarrow e'_U$ of  $T^1$-bundles such that
$f^*(\nabla_{e'_U})=\nabla_{e_U}$. The group of automorphisms of
$(e_U,\nabla_U)$ is the set of $T^1$-differentiable functions
defined on $U$.
 This
gerbe is trivial, since as remarked McDuff in [14], in this case
the Guillemin-Lerman-Sternberg method allows to construct a closed
form which extends $[\omega]$ if $H^1(F,{\R})=0$.

\bigskip

{\bf Remark}

 Suppose that the symplectic bundle $p:P\rightarrow N$ has a
Hamiltonian reduction. Then there exists an extension $\Omega$ of
$\omega$ (see [13]) which defines the distribution $D^{\Omega}$ on
$P$ as follows: let $u$ be an element of $P$,  $T_uP$ and
$TF_{p(u)}$  the respective tangent spaces of $P$ at $u$ and at
the fiber of $p(u)$.

$$
{D^{\Omega}}_u=\{v\in T_uP: \Omega_u(v,y)=0, y\in TF_{p(u)}\}.
 $$

 When the bundle is Hamiltonian, we can
suppose that the holonomy of the closed connection form is
Hamiltonian. And using a standard process, we can reduce the
structural group of this connection to its holonomy, and obtain
thus the Hamiltonian reduction.

\bigskip

{\bf Proposition 2.5.2.}

{\it Suppose that there exists an extension $[\Omega]$ of the
class $\omega$. Let $\Omega$ be a fixed representative. Then the
set of cohomology classes of closed $2$-forms $\Omega'$ whose
restriction to any fiber $F_u$ coincides with the restriction of
$\Omega$ to the fiber $F_u$ and such that $D^{\Omega}=D^{\Omega'}$
is isomorphic to $H^2(N,{\R})$.}

\bigskip

{\bf Proof.}

Remark that while this proposition is very similar to the problem
of the Brylinski's book [3] mentioned above, we cannot apply the
result obtained by Brylinski since we do not suppose that the
class $[\omega]$ is integral and $H^1(F,{\R})$ may not vanish.

Let $\Omega'$ be a representative of a cohomology class whose
restriction to the fiber $F_u$ of $p:P\rightarrow N$  coincide
with the restriction of $\Omega$ to $F_u$ and such that
$D^{\Omega}=D^{\Omega'}$. Then the form $\Omega-\Omega'$ projects
to a closed $2$-form $p(\Omega-\Omega')$ defined on the base, we
have thus defined a map between the set of extensions of
$[\omega]$ whose have a representant whose restriction to a fiber
$F_u$ coincides with the restriction of $\Omega$ to $F_u$ and
defines also the distribution $D^{\Omega}$, and $H^2(N,{\R})$ by
assigning to the class of $\Omega'$ the class of $\Omega-\Omega'$.
We have to show that this map is an isomorphism.

Suppose that the class of $p(\Omega-\Omega')$ is trivial. Then
there exists a $1$-form $\alpha$ on $N$ such that
$d(\alpha)=p(\Omega-\Omega')$. We denote  by $p^*(\alpha)$ the
pulls-back of $\alpha$ to $P$. This implies that
$\Omega'=\Omega+d(p^*(\alpha))$, thus the class of $\Omega$ and
$\Omega'$ coincide. This shows that the map $[\Omega']\rightarrow
[p(\Omega-\Omega')]$ is injective.

To show that this map is surjective, consider a closed $2$-form
$v$ of $N$, $p(\Omega-(\Omega-p^*(v)))=v$ $\bullet$

\bigskip

\medskip

The initial problem studied by Mc Duff was to find a Hamiltonian
reduction of the bundle $p:P\rightarrow N$, that is a symplectic
bundle  isomorphic to $p$,  whose transition functions take their
values in the Hamiltonian group of $(F,\omega)$. This problem can
be studied by a sheaf of categories. The definition of this sheaf
of category use the following result of Lalonde-McDuff [13]: which
allows to define its band:

\bigskip

{\bf Definition-Proposition 2.5.3.}

{\it Let $p:P\rightarrow N$ be a symplectic bundle. Suppose that
there exists a Hamiltonian reduction of $p$. Then there exists an
extension $\Omega$ of $\omega$, such that the Hamiltonian
reduction is defined by the holonomy of the closed connection form
$\Omega$. A Hamiltonian automorphism of the bundle $p:P\rightarrow
N$ is a diffeomorphism $\phi$ of $P$ which covers the identity,
such that the restriction of $\phi$ to the fiber over $n\in N$ is
an Hamiltonian automorphism of $(F,\omega_n)$, and such that
$\phi^*(\Omega)=\Omega$. We denote by $Aut(P,\Omega)$ the group of
Hamiltonian automorphisms of the Hamiltonian reduction
$(P,\Omega)$. The group $Aut(P,\Omega)$ does not depend of the
Hamiltonian reduction.}

\bigskip

{\bf Remark}

 In fact a more general
result is shown in Lalonde-McDuff [13] that is: the group of
diffeomorphisms $G(P,\omega)$ which cover the identity and such
that the restriction of each of its element $\phi$ to a fiber
$F_u$, belongs to the connected component of the group of
symplectic diffeomorphisms  of $(F_u,\omega_u)$, and which
preserves the symplectic class which defines the Hamiltonian
reduction does not depend of the chosen Hamiltonian reduction.
This result implies the one stated in the proposition above since
this group $G(P,\omega)$ contains $Aut(P,\Omega)$. The elements of
$Aut(P,\Omega)$ are the elements of $G(P,\Omega)$ which when
restricted to $(F_u,\omega_u)$ are Hamiltonian. We see that this
last condition is independent of the chosen Hamiltonian connection
$\Omega$ which defines any Hamiltonian reduction of
$p:P\rightarrow N$. A morphism $f:P\rightarrow P'$ between the
Hamiltonian bundles $P$ and $P'$ is a morphism of fiber bundles
$f$ such that $f^*(\Omega')=\Omega$, where $\Omega$ and $\Omega'$
are the closed connections forms whose holonomy define
respectively the Hamiltonian reduction of  $P$ and $P'$.

\bigskip

Now we can show the following:

\bigskip

{\bf Proposition 2.5.4.}

{\it Let $p:P\rightarrow N$ be a symplectic fibration. For any
open set $U$ of $N$, we define $C^1_F(\omega)(U)$ to be the
category whose objects are  Hamiltonian structures  on  the
symplectic bundle $p^{-1}(U)\rightarrow U$. A morphism between the
objects $(e_U,\Omega_U)$, and $(e',\Omega'_U)$ of
$C^1_F(\omega)(U)$  is a morphism of bundles $f:e_U\rightarrow
e'_U$ such that $f^*(\Omega'_U)=\Omega_U$. The correspondence
defined on the category of open subsets of $N$ by $U\rightarrow
C^1_F(\omega)(U)$ is a gerbe whose band $L$ is the sheaf induced
by the presheaf of Hamiltonian automorphisms such that for each
open set $U$ of $N$, and each $e_U$ of $C^1_F(\omega)(U)$, $L(U)$
is the group of Hamiltonian automorphisms of $e_U$ see
Definition-Proposition 2.5.3. The cohomology class of  the
classifying cocycle of $C^1_F(\omega)$ is the obstruction for the
existence of a Hamiltonian reduction of $p:P\rightarrow N.$}

\bigskip

{\bf Proof.}

\bigskip

Gluing conditions of objects.

\bigskip

 Consider  $(U_i)_{i\in I}$  an open cover of the open subset $U$
of $N$, such that $C^1_F(\omega)(U_i)$ is not empty, and
$(e_i,\Omega_i)$ an object of $C^1_F(\omega)(U_i)$. Suppose that
there exists a family of morphisms $u_{ij}:e^i_j\rightarrow e^j_i$
such that
${u_{i_1i_2}}^{i_3}{u_{i_2i_3}}^{i_1}={u_{i_1i_3}}^{i_2}$. Then
there exists a $F$-bundle $e$ over $U$ whose restriction to $U_i$
is $e_i$. We have to show that this bundle is Hamiltonian. Since
${u_{ij}}^*(\Omega_i)=\Omega_j$, the forms $\Omega_i$ glue
together to define on $e$ an extension $\Omega$ of $\omega$.
Consider a path $c:[0,1]\rightarrow N$, we can suppose that
$[0,1]$ is a union of intervals $I_l$ such that $I_l$ is contained
in $U_l$, an open set of the above cover. The holonomy of the
connection form $\Omega$ along $I_l$ is the holonomy of $\Omega_l$
along $I_l$, we conclude that the holonomy of $\Omega$ along $I$
is Hamiltonian, since each closed form $\Omega_l$ define an
Hamiltonian reduction on $e_l$.

\bigskip

Gluing conditions of arrows.

\bigskip

Let $e_U$ and $e'_U$ be a pair of objects of $C_F^1(\omega)(U)$.
The correspondence defined on the category of open subsets of $U$
which associates to $V$ the set of Hamiltonian morphisms
$Ham(e_U,e'_U)$ is a sheaf since it is the subsheaf of the sheaf
of morphisms between two bundles.

\bigskip

Connectivity.

 Let $e_U$ and $e'_U$ be a pair of objects of $C^1_F(\omega)$. We can suppose
that the open cover $(U_i)_{i\in I}$ of $U$ is a Hamiltonian
trivialization of the both bundles $e_U$ and $e'_U$. This implies
that the restrictions of $e_U$ and $e'_U$ to $U_i$ are isomorphic
as Hamiltonian bundles to the trivial Hamiltonian bundle
$U_i\times (F,\omega)$. We deduce that these Hamiltonian bundles
are locally isomorphic.

Let $(U_i)_{i\in I}$ be a symplectic trivialization  of
$p:P\rightarrow N$. The  trivial symplectic bundle $U_i\times
(F,\omega)$ is an element of $C^1_F(\omega)(U_i)$, which is not
empty.

The result of Lalonde and McDuff [13] recalled above shows that
the group $Aut(e_U,\Omega_U)$ of Hamiltonian automorphisms of the
Hamiltonian reduction of the restriction of $p$ to $p^{-1}(U)$
does not depend of the chosen object in $C^1_F(\omega)(U)$. This
implies that the correspondence defined on the category of open
subsets of $N$ by $U\rightarrow Aut(e_U,\Omega_U)$ defines a
presheaf $L'$ on $U$. We denote by $L$ the sheaf associated to
this presheaf. Remark that if $C^1_F(\omega)(U)$ is not empty,
then $L(U)=Aut(e_U,\Omega_U)$ for each object $e_U$ of
$C^1_F(\omega)(U)$. This implies that the gerbe is bounded by $L$
$\bullet$

\bigskip

{\bf 2.6 The McDuff construction of $Ham^s$, and closed connection
forms.}

\bigskip

The existence of a closed connection form $\Omega$ on the
symplectic bundle $p:P\rightarrow N$ does not insure the existence
of a Hamiltonian reduction of this bundle. This has motivated
McDuff to introduce the  group denoted $Ham^s$, such that the
existence of a closed connection form is equivalent to the
existence of a $Ham^s$-reduction. We will now present the
construction of the group $Ham^s$, and show using gerbe theory
that a $Ham^s$-reduction implies the existence of a closed
connection form on a symplectic bundle.

\bigskip

{\bf Definitions McDuff 2.6.1.}

Let $H_1(F,\omega,{\Z})$ be the first homology group of $F$ with
integral coefficients, we define $SH_1(F,\omega,{\Z})$ to be the
quotient of the integral $1$-cycles of $F$ by the image under the
boundary of $2$-cycles with zero symplectic area.  We denote
$SH_1(F,\omega,{\Q})$ to be the tensor product
$SH_1(F,\omega,{\Z})\otimes {\Q}$. Often we will respectively
denote $SH_1(F,\omega,{\Z})$, and $SH_1(F,\omega,{\Q})$ by
$SH_1(F,{\Z})$ and $SH_1(F,{\Q})$. Let $P_{\omega}$ be the values
of $\omega$ on rational cycles. We have the exact sequence:

$$
0\longrightarrow {\R}/P_{\omega}\longrightarrow
SH_1(F,{\Q})\longrightarrow H_1(F,{\Q})\longrightarrow 0
$$

Consider a section $s$ of $H_1(F,{\Q})\rightarrow SH_1(F,{\Q})$.
Then we can define on $Symp(F,\omega)$ the group of
symplectomorphisms of $(F,\omega)$, the map
$F_s:Symp(F,\omega)\rightarrow
H^1(F,{\R}/P_{\omega})=H^1(F,{\R})/H^1(F,P_{\omega})$ by:

$$
F_s(g)(u)=g(su)-s(gu)
$$
$\bullet$

 Recall that the group $Symp(F,\omega)$ acts canonically
on $SH_1(F,{\Q})$ and $H_1(F,{\Q})$. McDuff has shown that the
application $F_s$ is a $1$-cocycle for the canonical
representation defined on $Symp(F,\omega)$ which takes its values
in the group of linear automorphisms of
$H^1(F,{\R})/H^1(F,P_{\omega})$, and has defined $Ham^s$ to be the
kernel of this cocycle $F_s$.

\bigskip

{\bf Theorem McDuff 2.6.2.}

{\it A symplectic bundle $p:P\rightarrow N$ has a
$Ham^s$-reduction if and only if there exists a closed connection
form. Moreover, the group $Ham^s$ intersects every connected
component of $Symp(F,\omega)$.}

\bigskip

We will present now a proof of the first part of this theorem
using gerbe theory. In fact this problem can be reformulated in a
more general situation: Let $G$ be a Lie group whose dimension can
be infinite, and $H$ a subgroup of $G$, we suppose that $G/H$ is a
$K(\pi,1)$ space, that is its universal cover is contractible and
its fundamental group is $\pi$. We are looking for conditions
which insure the existence of a $H$-reduction. This problem can be
formulated using gerbe theory. We have:

\bigskip

{\bf Theorem 2.6.3.}

{\it Let $p:P\rightarrow N$ be a $G$-principal bundle defined on
$N$. Suppose either:

(i) the transitions functions $u_{ij}:U_i\cap U_j\rightarrow G$
take their values in the normalizer $Nor(H)$ of $H$ in $G$, where
$H$ is a subgroup of $G$, and $G/H$ is a $K(\pi,1)$ space.

(ii)or there exists a continuous representation $h:G\rightarrow L$
where $L$ is an abelian group isomorphic to the quotient of a
vector space $V$ by a discrete subgroup $\pi$  such that the
restriction of $h$ to the connected component  of the identity
$G_0$ of $G$ is trivial, a continuous surjective $1$-cocycle for
this representation whose kernel $H$ intersects every connected
component of $G$,

Then there exists a gerbe  $C_H$ defined on $N$, bounded by the
locally constant sheaf defined on $N$ by $\pi$ which represents
the obstruction of the bundle $p:P\rightarrow N$ to have a
$H$-reduction.}

\bigskip

{\bf Proof.}

The proof is a corollary of the following lemmas:

\bigskip

{\bf Lemma 2.6.4.}

{\it Suppose first that there exists a subgroup $H$ of $G$, such
that the transitions functions  $u_{ij}$ of $p:P\rightarrow N$ are
contained in the normalizer $Nor(H)$ of $H$ in $G$, then the right
quotient fiber by fiber of the bundle $p$ by $H$, defines a
$G/H$-bundle $\bar p:\bar P\rightarrow N$. Let $\hat {G/H}$ be the
universal cover of $G/H$. For each open subset $U$ of $N$, define
the category $C_H(U)$ to be the category whose objects are
$\hat{G/H}$-bundles whose quotient fiber by fiber by $\pi$ (recall
that $\pi$ is the fundamental group of $ {G/H}$) is the
restriction of $\bar p$ to $U$, a morphism $f:e_U\rightarrow e'_U$
between two objects $e_U$ and $e'_U$ of $C_H(U)$ is a morphism of
$\hat{G/H}$-bundles which projects to the identity on their
quotient by $\pi$. Then the correspondence defined on the category
of open subsets of $N$, $U\rightarrow C_H(U)$ defines a gerbe,
whose classifying cocycle is the obstruction for reduce the
structural group $G$ of the bundle $p$ to $H$.}

\bigskip

{\bf Proof.}

We have first to show the existence of the bundle $\bar p$. Let
$(U_i,u_{ij})$ be a trivialization of the bundle $p$. Since
$u_{ij}$ take their values in $Nor(H)$, for each element $x$ of
$U_i\cap U_j$, the right multiplication by $u_{ij}(x)$ of $G$
gives rise to a $G/H$-action of $u_{ij}(x)$ on $G/H$. We denote by
$\bar {u_{ij}}(x)$ this induced action. The map $\bar
{u_{ij}}:U_i\cap U_j\rightarrow G$, $x\longrightarrow
\bar{u_{ij}}(x)$ verified the Chasle relation, thus defines a
$G/H$-bundle $\bar p$ over $N$. Now we show that the
correspondence $U\rightarrow C_H(U)$ is a gerbe.

\bigskip

Gluing condition for objects.

\bigskip

Let $U$ be an open set of $N$,  $(U_i)_{i\in I}$ an open cover of
$U$, and $e_i$ an object of $C_H(U_i)$. We suppose that there
exists maps $g_{ij}:e^i_j\rightarrow e^j_i$ such that
${g_{i_1i_2}}^{i_3}{g_{i_2i_3}}^{i_1}={g_{i_1i_3}}^{i_2}$. Since
$e_i$ is a bundle, there exists a bundle $e$ over $U$ whose
restriction to $U_i$ is $e_i$. Since the restriction to $U_i$ of
the quotient fiber by fiber, of $e$ by $\pi$ is the quotient fiber
by fiber of $e_i$ by $\pi$, we deduce that $e$ is an element of
$C_H(U)$.

\bigskip

{Gluing condition for arrows.}

\bigskip

For each pair of objects $e$ and $e'$, of $C_H(U)$ the
correspondence defined on the category of open subsets of $U$ by
$V\rightarrow Hom(e_V,{e'}_V)$ where $e_V$ and ${e'}_V$ are the
respective restrictions of $e$ and $e'$ to $V$ defines a sheaf,
since it is the sheaf of morphisms between two bundles.

This shows that the correspondence defined on the category of open
subsets of $N$ by $U\rightarrow C_H(U)$ is a sheaf of categories.
Now we show that it is a gerbe.

\bigskip

Let $(U_i)_{i\in I}$ be a trivialization of the bundle $\bar p$,
then we can lift the restriction of $\bar p$ to $U_i$ to a bundle
$U_i\times \hat{G/H}$. This shows that $C_H(U_i)$ is not empty.

Let $U$ be an open set of $N$. Consider two objects $e_U$, and
$e'_U$ of $C_H(U)$. The restriction of $e_U$ and $e'_U$ to
$U_i\cap U$ are isomorphic to $U_i\cap U\times \hat{G/H}$ this
implies that the connectivity property holds.

The definition of $Hom(e_U,e_U)$, the group of automorphisms of an
object $e_U$ shows that its elements  coincide with the action of
 $\pi$, which thus defines a locally constant sheaf on $N$, which
is the band of $C_H$.

\bigskip

It remains to show that  the triviality of the gerbe $C_H$ is
equivalent to the existence of a $H$-reduction of $G$. Let $\hat
G$ and $\hat H$ be respectively the universal cover of $G$ and
$H$. The homotopy sequence applied to the fibration $\hat
H\rightarrow \hat G\rightarrow \hat G/\hat H$ implies that $\hat
G/\hat H$ is simply connected. The map $\hat G/\hat H\rightarrow
{G/H}$ is a covering map, thus $\hat G/\hat H$ is the universal
cover of $G/H$. Suppose that the gerbe  $C_H$ is trivial, then a
global object of this gerbe is a $\hat G/\hat H$-bundle. Since
$\hat G/\hat H$ is contractible, we deduce that this bundle is
trivial, and thus have a global section. This section projects to
a section of $\bar p$. This implies that the bundle $p$ has a
$H$-reduction $\bullet$

\bigskip

{\bf Lemma 2.6.5.}

{\it Suppose that there exists a continuous  representation
$h:G\rightarrow L$ (where $L$ is a quotient of a vector space $V$
by a discrete subgroup $\pi$), whose restriction to the connected
component $G_0$ of $G$ is trivial. Suppose also the existence of a
continuous cocycle $F$, surjective, for this representation whose
kernel $H$ intersects every connected component of $G$, then for
every principal $G$-bundle $p:P\rightarrow N$, there exists a
gerbe $C_H$, who represents the geometric obstruction for the
bundle $p$ to have a $H$-reduction.}

\bigskip

{\bf Proof.}

 For each elements $g\in G$, and $h\in H$, we have
$F(gh)=F(g)+gF(h)=F(g)$. This implies that the cocycle $F$ defines
a map $\bar F:G/H\rightarrow L$. The map $\bar F$ is surjective
since $F$ is surjective. Let $[g]$, and $[g']$ be two elements of
$G/H$, suppose that $\bar F([g])=\bar F([g'])$. Since $H$
intersects every connected component of $G$, we can choose two
elements $g$ and $g'$ in $G_0$, and respectively in the class
$[g]$ and $[g']$ such that $F(g)=F(g')$.

$$
F(g{g'}^{-1})=F(g)+h(g)F({g'}^{-1})=F(g)+F({g'}^{-1})=0
$$

since $g\in G_0$, and the restriction of $h$ to $G_0$ is trivial.
We deduce that $F(g)=F(g')$, thus $\bar F$ is a diffeomorphism.

Remark that $H\cap G_0=H_0$ is a normal subgroup of $G_0$. Since
$\bar F$ is a diffeomorphism, we deduce that $G/H=G_0/H_0$ is
diffeomorphic to $L$.

Consider now a $G$-bundle $p:P\rightarrow N$, defined by the
trivialization $u_{ij}:U_i\cap U_j\rightarrow G$. Then $F(u_{ij})$
defines a $L=G/H$-bundle $\bar p$ over $N$. The action of
$F(u_{ij}(x))$ on an element $[g]$ of $G/H$ is defined by
$F(gu_{ij}(x))$ where $g$ is an element of $[g]$ in $G_0$. For
each open subset $U$ of $N$, we define $C_H(U)$ to be the category
whose objects are $V$-bundles over $U$, whose quotient fiber by
fiber by $\pi$ is the restriction of $\bar p$ to $U$. Recall that
$L$ is the quotient of $V$ by $\pi$. The set of morphisms
$Hom(e,e')$ between two objects $e$ and $e'$ of $C_H(U)$ is the
set of morphisms of $L$-bundles which project to the identity on
the restriction of $\bar p$ to $U$. We are going to show that the
correspondence defined on the category of open subsets of $N$, by
$U\rightarrow C_H(U)$ is a gerbe who represents the geometric
obstruction to reduce $G$ to $H$.

\bigskip

Gluing property for objects.

\bigskip

Comsider an open subset $U$ of $N$, and an open covering
$(U_i)_{i\in I}$ of $U$. Let $e_i$ be an element of $C_H(U_i)$.
Consider a morphism $g_{ij}:e^i_j\rightarrow e^j_i$ such that
${g_{i_1i_2}}^{i_3}{g_{i_2i_3}}^{i_1}={g_{i_1i_3}}^{i_2}$. Since
$e_i$ are bundles, there exists a bundle $e$ over $U$ whose
restriction to $U_i$ is $e_i$. Since the restriction to $U_i$ of
the quotient fiber by fiber of $e$ by $\pi$ is the quotient fiber
by fiber of $e_i$ by $\pi$, we deduce that $e$ is an element of
$C_H(U)$.

\bigskip

{ Gluing condition for arrows.}

\bigskip

For each pair of objects $e$ and $e'$, the  correspondence defined
on the category of open subsets of $U$ by  $V\rightarrow
Hom(e_V,{e'}_V)$, where $e_V$ and ${e'}_V$ are the respective
restrictions of $e$ and $e'$ to $V$ defines a sheaf, since it is
the sheaf of morphisms between two bundles.

\bigskip

Let $(U_i)_{i\in I}$ be a trivialization of the bundle $\bar p$,
then we can lift the restriction of $\bar p$ to $U_i$, to a bundle
$U_i\times \hat{G/H}$. This shows that $C_H(U_i)$ is not empty.

Consider two objects $e_U$, and $e'_U$ of $C_H(U)$. The
restrictions of $e_U$ and $e'_U$ to $U_i\cap U$ are isomorphic to
$U_i\cap U\times \hat{G/H}$ this implies that the connectivity
property holds.

The definition of $Hom(e_U,e_U)$  the group of automorphisms of
the bundle $e_U$ shows that  it can be identified with $\pi$,
which thus defines a locally constant sheaf on $N$ which is the
band of $C_H$.

\bigskip

It remains to show that the triviality of the classifying cocycle
of the gerbe $c_H$ implies the existence of an $H$-reduction. Let
$\hat G$ and $\hat H$ be respectively the universal cover of $G$
and $H$. The homotopy sequence applied to the fibration $\hat
H\rightarrow \hat G\rightarrow \hat {G}/\hat H$ implies that $\hat
G/\hat H$ is simply connected. The map $\hat G/\hat H\rightarrow
\hat{G/H}$ is a covering map, thus $\hat G/\hat H$ is the
universal cover of $G/H$. Suppose that the gerbe  $C_H$ is
trivial, then a global object of this gerbe is a $\hat G/\hat
H$-bundle. Since $\hat G/\hat H$ is contractible, we deduce that
this bundle is trivial, and thus have a global section. This
section projects to a section of $\bar p$. This implies that the
bundle $p$ has a $H$-reduction $\bullet$

\bigskip

{\bf Remark.}

In the case of the lemma 2.6.5, above, the cohomology class of the
classifying cocycle $c_H$ is the obstruction for the bundle $\bar
p$ to be flat. This implies that it is the Chern class of this
bundle.

We are going to apply the above result to study the problem of the
existence of $Ham^s$-reductions.

\bigskip

{\bf Theorem 2.6.6.}

{\it Let  $p:P\rightarrow N$ be a symplectic bundle  whose typical
fiber is $(F,\omega)$,  then there exists a gerbe $C_{Ham^s}$ such
that the cohomology class $[c_{Ham^s}]\in H^2(N,
H^1(F,{\R})/H^1(F,P_{\omega}))$ of its classifying cocycle
$c_{Ham^s}$ is the obstruction to reduce the structural group of
the bundle to $Ham^s$. If the coordinate changes of the bundle
take their values in the connected component $Symp(F,\omega)_0$ of
$Symp(F,\omega)$, then there exists a gerbe $C_{Ham}$ whose
classifying cocycle is the obstruction for reduce the structural
group to $Ham(F,\omega)$.}

\bigskip

{\bf Proof.}

The group $Ham^s$ is the kernel of the continuous surjective
$1$-cocycle $F_s$, and it intersects every connected component of
$Symp(F,\omega)$. The quotient of $Symp(F,\omega)$ by $Ham^s$ is
$H^1(F,{\R})/H^1(F,P_{\omega})$. We can apply theorem 2.6.4.

Suppose that the coordinate changes take their values in
$Symp(F,\omega)_0$, since $Ham(F,\omega)$ is a normal subgroup of
$Symp(F,\omega)_0$, and the flux homomorphism allow us to identify
$Symp(F,\omega)_0/Ham(F,\omega)$ with
$H^1(F,\Gamma)/H^1(F,\Gamma)$, where $\Gamma$ is the flux group,
we can apply theorem 2.6.4 $\bullet$

\bigskip

{\bf Remark.}

In differential geometry, like in the theory of $G$-structures the
question of finding reductions of a $G$-bundle is intensively
studied. Let $H$ be a subgroup of $G$, if the left quotient $H/G$
is a $K(\pi,1)$ space, it is possible to write a similar theorem
to the one above, and obtain an obstruction cocycle whose
cohomology class decide of the existence of a $H$-reduction. This
can be for example  applied to the existence of a riemannian
structure on a manifold, and also to solve differential equations
defined on jet-bundles, since in many cases the existence of
solutions is equivalent to the existence of reductions of jets
bundles.

\bigskip

We will give now another proof of the theorem of McDuff mentioned
above  which says that the existence of a closed connection form
on a symplectic bundle $p:P\rightarrow N$ implies the existence of
a $Ham^s$-reduction.

\bigskip

{\bf Theorem McDuff 2.6.7. (see [14]).}

{\it Let $p:P\rightarrow N$ be a symplectic bundle endowed with a
closed connection form, then there exists on $P$ a
$Ham^s$-reduction.}

\bigskip

{\bf Other proof.}

Suppose the existence of a closed connection form defined on the
bundle $p:P\rightarrow N$. We have to show that the cohomology
class $[c_{Ham^s}]$ is trivial. It has been shown by
McDuff-Lalonde [13], that the holonomy around a contractible loop
is Hamiltonian. Consider the reduction of the symplectic bundle to
the holonomy of the closed connection form. Since the Hamiltonian
group is the connected component of $Ham^s$, we deduce that the
composition of the transitions functions $u_{ij}$ and of $F_s$,
$F_s(u_{ij})$ is constant, if needed, we  shrink the open set
$U_i$ such that $u_{ij}(U_i\cap U_j)$ is contained in the same
connected component of $Symp(F,\omega)$. This implies that the
bundle $\bar p$ (defined in the proof of Lemma 2.6.5) is flat.
Thus its Chern class is a torsion class. Since the lattice $\pi$
in this case is a ${\Q}$-vector space, we deduce that the Chern
class of this bundle is zero $\bullet$

\bigskip

{\bf Sketch of the proof of McDuff [14].}

McDuff defines for each symplectic bundle $p:P\rightarrow N$ of
fiber $(F,\omega)$, a cohomology $2$-class in
$H^2(N,H^1(F,P_{\omega}))$ (in fact it is the class of 2.6.5) as
follows: The bundle $p$ is defined by a classifying map
$p':N\rightarrow BSymp(F,\omega)$. The map $F_s$ induces a map
$F'_s:BSymp(F,\omega)\rightarrow BH^1(F,{\R})/H^1(F,P_{\omega})$.
There exists a $Ham^s$-reduction if and only if the composition
$F'_s\circ p'$ is null homotopic, since we have an exact sequence

$$
1\longrightarrow Ham^s\longrightarrow
Symp(F,\omega)\longrightarrow
H^1(F,{\R})/H^1(F,P_{\omega})\longrightarrow 1.
$$

The space $BH^1(F,{\R})/H^1(F,P_{\omega})$ is a
$K(H^1(F,P_{\omega}),2)$-space, and the set of homotopy classes of
maps $N\rightarrow K(H^1(F,P_{\omega}),2)$ is one to one with
$H^2(N,P_{\omega})$. The obstruction class of McDuff is defined to
be the homotopy class of $F'_s\circ p'$.

 The proof of McDuff of the previous result, is done by showing that the
previous class vanishes on the $2$-sub-complex of the
$CW$-complexe $N$. On this purpose, she shows that it is the image
by a null connecting homomorphism related to an exact sequence of
a one class $\bullet$

\bigskip

{\bf 2.7 The universal obstruction of McDuff.}

\bigskip

In this part, we will show how the universal class defined by
McDuff can be defined using gerbe theory.

Let $ESymp(F,\omega)\rightarrow BSymp(F,\omega)$ be the universal
bundle of the group $Symp(F,\omega)$.   The $1$-cocycle
$F_s:Symp(F,\omega)\rightarrow H^1(F,{\R})/H^1(F,P_{\omega})$
defined by McDuff, induces a
$H^1(F,{\R})/H^1(F,P_{\omega})$-bundle on $BSymp(F,\omega)$. See
lemma 2.6.5. The Chern class $U_F$ of this bundle is the universal
class $U_M$ defined by McDuff, it can be viewed as the cohomology
class of the classifying cocycle of the gerbe which represents
geometrically the obstruction for the previous
$H^1(F,{\R})/H^1(F,P_{\omega})$-bundle to be trivial. Since each
$(F,\omega)$-symplectic bundle $p:P\rightarrow N$, is classified
by a classifying map $f:N\rightarrow BSymp(F,\omega)$, the
obstruction class to obtain a $Ham^s$-reduction is $f^*(U_F)$.
This is the class defined in the sketch of the proof of McDuff in
2.6.

\bigskip

{\bf 2.8 Generalizations to topoi.}

\bigskip

The previous construction applied to symplectic bundles can be
generalized to other situations; algebraic geometry,
arithmetic,... On this purpose we will adapt this result to topoi.

\bigskip

{\bf Definitions 2.8.1.}

Let $G$ be a group endowed with a topology, the topology can be
the Zariski, etale, the Lie topology, ect... A continuous right
$G$-action of $G$ on the topos $(P,J_P)$, is a continuous functor
$d_G:P\times G\rightarrow P$, such that if $u$ is the
multiplication of $G$ by $g$, $d_G\circ (Id_P\times
u)=d_G(d_G\times Id_G)$.

A {\bf $G$-torsor} defined on a topos $N$ is a continuous functor
$p:(P,J_P)\rightarrow (N,J_N)$ such that:

 (i) $(P,J_P)$ is endowed with an action of $G$, $p$ commutes with the
action of $G$ that is, the composition $P\times G \rightarrow
P\rightarrow N$, (where $P\times G\rightarrow P$ is the canonical
projection, and $P\rightarrow N$ is $p$) and $P\times G\rightarrow
P\rightarrow N$ (where $P\times G\rightarrow P$ is the
multiplication $d_G$ and $P\rightarrow N$ is $p$) coincide.

 (ii) The canonical map $P\times G\rightarrow P\times P\times
G\rightarrow  P\times P$ which is the composition of the canonical
embedding $P\times G\rightarrow P\times P\times G$, and the
product of the identity on the first factor, and the
multiplication $d_G$ on the second and third factor is an
isomorphism. We suppose that the quotient of $P$ by $G$ is $N$.
Recall that the quotient of $P$ by $G$ is the initial element in
the category of maps $p':P\rightarrow N'$ such that $p'$ commutes
with the action of $G$ $\bullet$

We will assume that the torsor is locally trivial. This means that
there exists a covering family of $N$, $(U_i)_{i\in I}$ such that:

There exists an isomorphism $u_i:P_{\mid U_i}\rightarrow U_i\times
G$ between the restriction $P_{\mid U_i}$ of $P$ to $U_i$ and
$U_i\times G$. We can thus define $u'_{ij}={u_i\circ
{u_j}^{-1}}_{\mid U_i\times_NU_j\times G}:U_i\times_NU_j\times
G\rightarrow U_i\times_NU_j\times G$. Let $e':G\rightarrow G$,
$g\rightarrow e$, where $e$ is the neutral of $G$ and
$e_{ij}:U_i\times_NU_j\times G\rightarrow G$ the canonical
projection. We can define $u_{ij}:U_i\times_NU_j\rightarrow G$ by
$e_{ij}\circ u'_{ij}\circ (Id_{U_i\times_NU_j}\times e')$. We have
 ${u'_{i_1i_2}}^{i_3}{u'_{i_2i_3}}^{i_1}={u'_{i_1i_3}}^{i_2}$; $P$
is obtained by gluing  the family of $(U_i\times G)_{i\in I}$
using $u'_{ij}$ .

Let $H$ be a subgroup of $G$, we say that the torsor $P\rightarrow
N$ has a $H$-reduction if and only if it is isomorphic to a torsor
whose transition functions $u_{ij}$ take their values in $H$.

\bigskip

{\bf Theorem 2.8.2.}

{\it  Let $p:P\rightarrow N$ be a $G$-torsor. Suppose that either,

1. there exists a subgroup $H$ of $G$ such that, $G/H$ is a
$K(\pi,1)$ space, and the torsor has a  $Nor(H)$-reduction.

2. or  there exists a $1$-cocycle surjective and continue
$F:G\rightarrow L$ for a representation $h$ of $G$, where $L$ is
the quotient of a vector space by a discrete subgroup $\pi$, such
that the restriction of $h$ to the connected component $G_0$ of
$G$ is trivial, and the kernel $H$ of $F$ intersects every
connected component of $G$,

3. or the left quotient $H/G$ is a $K(\pi,1)$-space,

then there exists a gerbe $C_H$ defined on $N$ such that the
cohomology class of its classifying cocycle is the obstruction for
reduce $G$ to $H$.}

\bigskip

{\bf Proof.}

We will only give the proof  in the first case. The fact that
 the torsor has a $Nor(H)$-reduction implies the existence of a
 $G/H$-torsor $\bar P$, which is the right quotient of $P$ by $H$.

For each object $U$ of $N$, we define $C_H(U)$ to be the category
 whose objects are $\hat {G/H}$-torsors whose quotient by $\pi$
is the restriction of $\bar P$ to $U$. A morphism between two
objects of $C_H(U)$, is a morphism of torsors which projects to
the identity on the restriction of $\bar P$ to $U$.

Now we show that the correspondence $U\rightarrow C_H(U)$ is a
gerbe.

\bigskip

Gluing condition for objects.

\bigskip

Let $U$ be an object of $N$,  $(U_i)_{i\in I}$ a covering family
of $U$, $e_i$ an object of $C_H(U_i)$. We suppose that there
exists maps $g_{ij}:e^i_j\rightarrow e^j_i$ such that
${g_{i_1i_2}}^{i_3}{g_{i_2i_3}}^{i_1}={g_{i_1i_3}}^{i_2}$. Since
$e_i$ are torsors, there exists a torsor $e$ over $U$ whose
restriction to $U_i$ is $e_i$. Since the restriction to $U_i$ of
the quotient of $e$ by $\pi$ is the quotient  of $e_i$ by $\pi$,
we deduce that $e$ is an element of $C_H(U)$.

\bigskip

{ Gluing condition for arrows.}

\bigskip

For each objects $e$ and $e'$, the set of morphims defined on the
sub-topos over   $U$, by  $V\rightarrow Hom(e_V,{e'}_V)$, where
$e_V$ and ${e'}_V$ are the respective restrictions of $e$ and $e'$
to $V$ defines a sheaf of sets, since it is the sheaf of morphisms
between two torsors.

\bigskip

Let $(U_i)_{i\in I}$ be a trivialization of the torsor $\bar P$,
we can lift the restriction of $\bar P$ to $U_i$ to the torsor
$U_i\times \hat{G/H}$. This shows that $C_H(U_i)$ is not empty.

Consider two objects $e_U$, and $e'_U$ of $C_H(U)$. The
restrictions of $e_U$ and $e'_U$ to $U_i\times_N U$ are isomorphic
to $U_i\times_N U\times \hat{G/H}$ this implies that the
connectivity property holds.

The group $Hom(e_U,e_U)$ is the the group of automorphisms of the
torsor $e_U$ which project to the identity isomorphism of the
restriction of $\bar p$ to $U$. This group is identified to $\pi$,
which thus defines a locally constant sheaf on $N$ which is the
band of $C_H$ $\bullet$

\bigskip

{\bf Remark.}

The triviality of the gerbe $C_H$ does not necessarily implies the
existence of a $H$-reduction, if $N$ is not a manifold. Since for
other categories,  homotopy is not well-understood, there are no
precise definitions of null-homotopic maps.

\bigskip

{\bf 3. The group $Ham^s$ and the etale topos of a manifold.}

\bigskip

The group $Ham^s$ introduced by McDuff allows to characterize
symplectic bundles whose have a closed connection form  to be the
symplectic bundles  endowed with a $Ham^s$-reduction. In [14] it
is shown that a $Symp(F,\omega)_0$-bundle $p:P\rightarrow N$ is
endowed with a closed connection form if and only if  there exists
a finite cover $\hat N$ of $N$, such that the pull-back of $p$ to
$\hat N$ has a Hamiltonian reduction. This motivates to define
$Symp(F,\omega)$-bundles on the etale topos of $N$. The motivation
is due to this historical remark: In algebraic geometry, algebraic
principal bundles are locally trivial up to a finite etale cover.
This has motivated the definition of the etale topology.

 \bigskip

 {\bf Definitions 3.1.}

 The {\bf etale topos} of a manifold $N$ is the category whose objects
 are differentiable maps $c:U\rightarrow N$  which
 are finite covering maps  onto their images. A morphism between two
 objects is a  covering map.

A covering family of the etale topos, $Et(N)$ of $N$, is a family
$(U_i)_{i\in I}$ such that the arrow $u_i:U_i\rightarrow N$ is a
finite etale cover, and the union of $(u_i(U_i))_{i\in I}$ is $N$

 A symplectic bundle  $p:P\rightarrow Et(N)$ whose typical fiber is the symplectic
manifold $(F,\omega)$ is defined by a covering family $(U_i)_{i\in
I}$ of $Et(N)$ for the etale topology. The transition functions
are symplectic bundles isomorphisms  of the trivial symplectic
bundle $U_i\times_NU_j\times Symp(F,\omega)$, defined by
$u_{ij}:U_i\times_NU_j\rightarrow Symp(F,\omega)$ such that
${u_{i_2i_3}}^{i_1}{u_{i_1i_2}}^{i_3}={u_{i_1i_3}}^{i_2}$.

A closed connection form on the symplectic bundle is  defined by a
family of closed connections forms $\Omega_i$  of the bundle
$e_i:U_i\times(F,\omega)$ (recall that $\Omega_i$ is a $2$-form
which extends $\omega$), such that on $U_i\times_NU_j$, we have:
${u_{ij}}^*({\Omega_i}_{\mid U_i\times_NU_j})={\Omega_j}_{\mid
U_i\times_NU_j}$, where ${\Omega_i}_{\mid U_i\times_NU_j}$ and
${\Omega_j}_{\mid U_i\times_NU_j}$ are the respective restrictions
of $\Omega_i$ and $\Omega_j$ to $U_i\times_NU_j$.

A symplectic bundle defined on $N$ induces canonically a
symplectic bundle on $Et(N)$, since an open covering of $N$
defines an etale covering of $N\bullet$

\bigskip

{\bf Proposition 3.2.}

{\it Let $P$ be a symplectic bundle defined on the etale topos of
a manifold $N$, then there exists a symplectic bundle $\hat P$
defined on a covering space $\hat N$ of $N$, such that  the
symplectic bundle induced by $\hat P$ on $Et(\hat N)$, is the
pull-back of $P$ by the covering map $\hat N\rightarrow N$. If $N$
is compact we can suppose that $\hat N$ is a finite cover.}

\bigskip

{\bf Proof.}

Let $(d_i:U_i\rightarrow N)_{i\in I}$ be the etale covering family
of $N$ which defines the symplectic bundle. Then we can define a
manifold $\hat N$ as follows: $\hat N$ is the quotient of the
union of $U_i$ by identifying the elements $u_i\in U_i$, and
$u_j\in U_j$ such that $d_i(u_i)=d_j(u_j)$. We denote by
$l_i:U_i\rightarrow \hat N$ the canonical map. The manifold $\hat
N$ is a cover of $N$ since the restriction of the canonical
projection $\hat N\rightarrow N$ to $l_i(U_i)$  is
$d_i{l_i}^{-1}$.

There exists a diffeomorphism $l_{ij}:l_i(U_i)\cap
l_j(U_j)\rightarrow U_i\times_NU_j$, such that on
$l_{i_1}(U_{i_1})\cap l_{i_2}(U_{i_2})\cap l_{i_3}(U_{i_3})$,
${{l_{i_2i_3}}^{i_1}}^{-1}{l_{i_1i_2}}^{i_3}=Id_{l_{i_1}(U_{i_1})\cap
l_{i_2}(U_{i_2})\cap l_{i_3}(U_{i_3})}$, thus we can define the
symplectic bundle $\hat P$ on $\hat N$ by gluing $l_i(U_i)\times
(F,\omega)$ using ${u'}_{ij}= u_{ij}\circ l_{ij}$, where $u_{ij}$
are the transition functions of $P$.

The construction of $\hat P$ shows that the induced bundle on
$Et(\hat N)$, by $\hat P$, is the pull-back of $P$ by the
canonical map $\hat N\rightarrow N$. If $N$ is compact, then we
can suppose that there exists a finite number of $U_i$, this
implies that $\hat N$ is compact, therefore is a finite cover of
$N$ $\bullet$

\bigskip

We can rewrite the theorem of McDuff [14] as follows:

\bigskip

{\bf Theorem 3.3.}

{\it Let $p:P\rightarrow Et(N)$ be a $Sym(F,\omega)_0$-bundle
defined on the etale topos of a compact manifold $N$, then $P$ has
a closed connection form if and only if it has a Hamiltonian
reduction.}

\bigskip

{\bf Proof.}

The previous proposition shows that there exists a finite cover
$\hat N$ of $N$, and an induced symplectic bundle $\hat P$ over
$\hat N$. Suppose that the closed connection form is defined on
$P$ by the family of $2$-forms $\Omega_i$  defined on the etale
cover $(l_i:U_i\rightarrow N)_{i\in I}$. As in the previous
proposition, we can show that there exists a finite cover $N'$ of
$N$ such that the pull-back $P'$  of $\hat P$ to $N'$ is endowed
with a closed connection form, such that the closed connection
form induced on its etale cover is defined on $U'_i=P'\times_NU_i$
by the pull-back of $\Omega_i$ by $P'\times_NU_i\rightarrow U_i$.
We can apply the result of McDuff and obtain a Hamiltonian
reduction
 $P"\rightarrow N"$ on the pull-back of $P"$ of $P'$ to a finite cover $N"$ of $N'$.
  We denote by $l"_i$ the canonical map
  $l"_i:U"_i=U'_i\times_{N'}P"\rightarrow N"$.
There exists a family of maps $u"_i:l"_i(U"_i)\rightarrow
Symp(F,\omega)$ such that $u"_iu"_{ij}{u"_j}^{-1}\in
Ham(F,\omega)$, where $u"_{ij}$ are  the coordinate changes of $
P"$, thus $u"_iu"_{ij}{u"_j}^{-1}{l"_j}^{-1}$ defined a
Hamiltonian reduction of $P$. Since the family $(U"_i\rightarrow
N")_{i\in I}$ is an etale cover of $N"$, $(U"_i\rightarrow
N"\rightarrow N)_{i\in I}$ is also an etale cover of $N$ $\bullet$

\bigskip

{\bf 4. Flux, and  holonomy of gerbes.}

\bigskip

In this part, we will relate the flux of a symplectic manifold
$(F,\omega)$ to the  holonomy of the gerbe $C(\omega)$ defined in
2.4.

\bigskip

Let $E$ be a $T^1$-gerbe defined on $P$, that is a gerbe such that
for each open set $U$ of $P$, $E(U)$ is a category of
$T^1$-bundles defined on $U$. Consider an open covering
$(U_i)_{i\in I}$ of $P$, such that $U_i$ is contractible. Let
$e^j_i$ be the restriction of an object $e_i$ of $E(U_i)$ to
$U_i\cap U_j$, there exists a morphism $u_{ij}:e^i_j\rightarrow
e^j_i$. We denote by $c_{ijl}$, the automorphism
$u_{li}u_{ij}u_{jl}$ of the restriction of  $e_l$ to $U_i\cap
U_j\cap U_l$. It is defined by a $T^1$-differentiable function.
Since $c_{ijl}$ is the classifying $2$-cocycle of $E$, there
exists a $1$-chain $h_{ij}$ of $1$-forms such that:

$$
h_{jl}-h_{il}+h_{ij}=-{i\over{2\pi}}d(Log(c_{ijl}))
$$

since $d(h_{ij})$ is a $1$-cocycle, there exists a $0$-chain of
$2$-forms $L_i$ such that

$$
L_j-L_i=d(h_{ij})
$$

\bigskip

{\bf Definitions 4.1.}

The family of forms $h_{ij}$ is a called a {\bf connection of the
gerbe}, and the family of forms $(L_i)_{i\in I}$ is the curving of
the gerbe, this means that
 there exists a related connective structure $Co$
defined on the gerbe,  and elements $\alpha_i$ of $Co(e_i)$, such
that $h_{ij}=\alpha_j-{u_{ij}}^*\alpha_i$.
 The $3$-form
whose restriction to $U_i$ is $dL_i$ is the curvature of the
connective structure.
 Suppose that the curvature is zero, then
$L_i=d(L'_i)$, $h_{ij}=L'_j-L'_i+d(h'_{ij})$, we denote by
$c'_{i_1i_2i_3}=-{i\over
{2\pi}}Log(c_{i_1i_2i_3}^{-1})+{h'}_{i_2i_3}-{{h'}_{i_1i_3}}+{h'}_{i_1i_2}$
to be {\bf the holonomy of the connection}, $c'_{i_1i_2i_3}$  is
constant and is a $2$-cocycle $\bullet$

\bigskip

{\bf Definition 4.2.}

 For each map $l:N_2\rightarrow P$, where $N_2$ is a surface without a
 boundary, the pull-back of  the gerbe, and its connective structure to $N_2$, by $l$ has a
 vanishing curving.
 Using the Cech-De Rham isomorphism, we can identify the holonomy
 cocycle of this gerbe
 with a $2$-form $Hol(h_{ij},N_2)$. The holonomy of the connection on $N_2$ is

 $$
 \int_{N_2}Hol(h_{ij},N_2)
 $$

$\bullet$
\bigskip

 Let $(F,\omega)$ be a symplectic manifold, and $C_F(\omega)$ the
$T^1$-gerbe representing the obstruction of $[\omega]$ to be
integral. If the   band of this gerbe is extended to the sheaf of
differentiable $T^1$-functions, it becomes  trivial and flat.

 For each open set $U$ of $F$, the set of connections defined on an
object $e_U$ of $C_F(\omega)(U)$ which curvature is the
restriction of $\omega$ to $U$ defines a connective structure, the
curving of the connective structure is the restriction of $\omega$
to $U_i$. The cocycle representing the holonomy of this connective
structure is the image of $\omega$ by the Cech-De Rham
isomorphism. This can be deduced from 2.4.3.

\bigskip
\medskip

Let $l:N_2\rightarrow F$ be a differentiable map defined on the
surface $N_2$, the holonomy of this connective structure around
$N_2$ is:

$$
\int_{N_2}l^*(\omega)
$$

 This
definition is related to the definition of the flux, since for
each path $\gamma=c_t$ of $N$, and each path $\phi_t$ of the
connected component of $Symp(F,\omega)$, $\phi_t(\gamma)$ is a map
from $h:I_2\rightarrow F$, the flux of $\phi_t(\gamma)$ is nothing
but the half of the holonomy around the sphere $S^2$ obtained by
gluing two copies of $I_2$ along their boundaries. The map
$f:S^2\rightarrow N$ is obtained by restricting $h$ to each copy
of $I_2$. The holonomy of $f$ is defined to be the limit of the
holonomy of a sequence of differentiable maps which converges
towards $f$.

\bigskip

{\bf 5. A geometric interpretation of a section of
$H_1(M,{\R})\rightarrow SH_1(M,{\R})$.}

\bigskip

 In [14], McDuff gives  a geometric
interpretation of a section $p:H_1(M,{\R})\rightarrow
SH_1(M,{\R}),$ when the cohomology class $[\omega]$ is integral.
In this section we generalize this interpretation when $[\omega]$
is not necessarily integral. We denote by
$\pi:SH_1(M,{\R})\rightarrow H_1(M,{\R})$ the projection map.
Suppose that the class $[\omega]$, of the symplectic manifold
$(F,\omega)$ is not necessarily integral. Consider a cycle
$[\gamma]$ represented by the chain $h:T^1\rightarrow F$, where
$T^1$ is the circle, the pull-back by $h$, of the gerbe
$C_F(\omega)$, to $T^1$ is trivial.

\bigskip

{\bf Proposition 5.1.}

{\it  Consider an object $e$ of $h^*(C_F(\omega))$ which is the
pull-back of an object $e'$ of a tubular neighborhood of $h(T^1)$.
Let $L$ be a connection in  $ h^*Co(e')$. Denote by $h'_L(\gamma)$
the holonomy around $\gamma$ of $L$. It does not depend of the
element chosen in $h^*(Co(e'))$.}

\bigskip

{\bf Proof.}

 To show this, consider
another connection $L'$ in $h^*(Co(e'))$. We can suppose that
$h(T^1)$ is covered by $(U_i)_{i\in I}$,  the union of $U_i$ is a
tubular neighborhood of $h(T^1)$, and $C_F(\omega)(U_i)$ is not
empty. The fact that the union of $U_i$ is a tubular neighborhood
of $h(T_1)$ implies that
 the restrictions of $L$ and $L'$ to $I_i=h^{-1}(h(T^1)\cap U_i)$
can be supposed to be  the pull-back of elements $d+\alpha_i$ and
$d+\alpha'_i$ of $Co(e_i)$, where $d$ is the differential  $e_i$
is an object of $C_F(\omega)(U_i)$ and $\alpha_i$ and $\alpha'_i$
are $1$-forms defined on $U_i$. We have $\alpha'_i=\alpha_i
+df'_i$ where $f'_i$ is a function defined on $U_i$ since the
curving of the gerbe is the closed form $\omega$. Denote by $L_i$
and $L'_i$ the restrictions of $L$ and $L'$ to $I_i$. On $I_i$,
$L_i=d+du_i$, the coordinate changes $v_{ij}$ of the bundle $e$
are defined by $du_j-du_i=-i{1\over {2\pi}}dLog(v_{ij})$, the
holonomy cocycle of $L$ is given by
$-{i\over{2\pi}}Log(v_{ij}^{-1})-u_j+u_i)$. Since
$L'_i=d+d(u_i+f_i)$, where $f_i$ is the pulls-back of $f'_i$ by
$h$. We deduce that the holonomy cocycle of $L$ and $L'$ coincide
up to a boundary. Thus their cocycle have the same cohomology
class $\bullet$

\bigskip

 We can define
$h_L(\gamma)$ the image of the holonomy of this connection in
${\R}/P_{\omega}$. Let $[\gamma]\in H_1(F,{\R})$,  defines the
section $p([\gamma])$ to be the class of elements $\gamma'$ in
$\pi^{-1}([\gamma])$ such that the holonomy  around $\gamma'$ is
in $ P_{\omega}$.

\bigskip

 {\bf 6. Existence of symplectic bundles and gerbes.}

\bigskip

Let $F$ be the flux, and $\Gamma_{\omega}$ be the flux group. The
flux conjecture has been shown recently by Ono, thus
$\Gamma_{\omega}$ is a discrete subgroup of $H^1(F,{\R})$. There
exists an exact sequence

$$
1\longrightarrow Ham(F,\omega)\longrightarrow
Symp_0(F,\omega)\longrightarrow
H^1(F,{\R})/\Gamma_{\omega}\longrightarrow 1.
$$

Let $p:P\rightarrow N$ be a symplectic bundle defined by the
cordinate changes $g_{ij}:U_i\cap U_j\rightarrow Symp_0(F,\omega)$
on the trivialization $(U_i)_{i\in I}$. We can project the cocycle
$g_{ij}$ to maps $F(g_{ij})={g'}_{ij}:U_i\cap U_j\rightarrow
H^1(F,{\R})/\Gamma_{\omega}$, and obtain a
$H^1(F,{\R})/\Gamma_{\omega}$-bundle as in 2.6.5. A natural
question is the following: given a
$H^1(F,{\R})/\Gamma_{\omega}$-bundle $\bar p$ is there a
symplectic bundle which gives rise to $\bar p$.

\bigskip

This problem is an example of the basic examples which have
motivated the definition of gerbes theory. Consider an open set
$U$ of $N$, and $C(U)$ a category of symplectic bundles whose
transition functions take their values in $Symp(F,\omega)_0$, and
which induces the restriction of $\bar p$ to $U$. Suppose that
$\bar p$ is defined by the transition functions $g'_{ij}$, and
there exists elements $g_{ij}$ over $g'_{ij}$ such that the
conjugation by $g_{ij}$ in $Ham(F,\omega)$ defined a bundle over
$N$ whose typical fiber is $Ham(F,\omega)$. We suppose also that
the automorphisms group of an object $e_U$ of $C(U)$ are the
sections of the previous $Ham(F,\omega)$-bundle. we denote $L_1$
the sheaf of those sections. The correspondence $U\rightarrow
C(U)$ is a gerbe bounded by $L_1$.

\bigskip

Denote by $l$ the rank of the group $\Gamma_{\omega}$, the torus
$T^l$ is the maximal compact subgroup of
$H^1(M,{\R})/\Gamma_{\omega}$. The bundles defined over $N$, which
fiber is $T^l$ are classified by their first Chern class. This can
enable to construct symplectic bundles which does not admit
Hamiltonian reductions if the Chern class is not zero.

\bigskip

{\bf 7. $2$-gerbes, $2$-gerbed towers.}

\bigskip

The notion of $2$-gerbe has been defined by Lawrence Breen [2],
[3] it allows to represent geometrically $3$-cohomology classes.
In the preprint [20], Tsemo has defined the notion of gerbed
towers, this is a recursive definition of geometric
representations of cohomology classes. We will present now the
notion of $2$-gerbes, and $2$-gerbed towers, which enable us to
cope with the extension problem when $[\omega]$ is not necessarily
an integral class. An alternative discussion has been presented
above using the group $Ham^s$, the construction given in this
section allows to show the existence of a connection on a bundle
which has a Hamiltonian reduction without using the
Guillemin-Lerman-Sternberg construction. The definition  of sheaf
of $2$-categories uses the definition of $2$-categories or
bicategories which has been defined by Benabou.

\bigskip

{\bf Definition.}

\bigskip

A {\bf bicategory} $C$ is defined by a class of objects $C$, for
each pair of objects $u$ and $v$ of $C$, a category  $Hom(u,v)$.
The objects of $Hom(u,v)$ are called the $1$-arrows, and the
arrows of $Hom(u,v)$ are the $2$-arrows, there exists a
composition map:

$$
Hom(u_2,u_3)\times Hom(u_1,u_2)\longrightarrow Hom(u_1,u_3)
$$

For each quadruple $(u_1,u_2,u_3,u_4)$, there exists an
isomorphism between the functors

$$
(Hom(u_3,u_4)\times Hom(u_2,u_3))\times
Hom(u_1,u_2)\longrightarrow Hom(u_1,u_4)
$$

and

$$
Hom(u_3,u_4)\times (Hom(u_2,u_3)\times
Hom(u_1,u_2))\longrightarrow Hom(u_1,u_4)
$$

which satisfies more compatibility axioms which can be found in
Breen [2] $\bullet$

\bigskip

{\bf Definitions.}

 Let $N$ be a manifold, a
{\bf sheaf of $2$-categories} is a correspondence $C$ defined on
the category of open subsets of $N$ by:

$$
U\longrightarrow C(U)
$$

where $C(U)$ is a $2$-category, which verifies the following
properties:

for each embedding map $U\rightarrow V$, there exists a
restriction functor $r_{U,V}:C(V)\longrightarrow C(U)$, such that

$$
r_{U_1,U_2}\circ r_{U_2,U_3}=r_{U_1,U_3}.
$$

\bigskip

{\bf Gluing properties for objects.}

\bigskip

Let $(U_i)_{i\in I}$ be a covering family of  an open set $U$ of
$N$,  $e_i$ an object of $C(U_i)$, and a $1$-arrow
$g_{ij}:r_{U_i\cap U_j,U_j}(e_j)\rightarrow r_{U_i\cap
U_j,U_i}(e_i)$, suppose there exists a $2$-arrow
$h_{i_1i_2i_3}:{g_{i_1i_2}}^{i_3}{g_{i_2i_3}}^{i_1}\rightarrow
{g_{i_1i_3}}^{i_2}$ which satisfies:

$$
{h_{i_1i_2i_4}}^{i_3}(Id\circ
{h_{i_2i_3i_4}}^{i_1})={h_{i_1i_3i_4}}^{i_2}({
h_{i_1i_2i_3}}^{i_4}\circ Id)
$$

then there exists an object $e$ of $C(U)$ whose restriction to
$U_i$ is $e_i$.

\bigskip

{\bf Gluing conditions for arrows.}

\medskip

For each pair of objects $e$ and $e'$ of $U$, the correspondence
defined on the category of open sets contained in $U$ by
$V\rightarrow Hom(r_{U,V}(e),r_{U,V}(e'))$ defines a sheaf of
categories.

\bigskip

A {\bf $2$-gerbe} is a sheaf of bicategories which satisfies the
following:

 1. The bicategory $C(U)$ is a $2$-groupoid, this means that
$1$-arrows are inverible up to $2$-arrows, and $2$-arrows are
invertible.

 2. For every point $x$ of $N$, there exists a neighborhood $U_x$ of $x$,
such that $C(U_x)$ is not empty.

 3. Any pair of objects $e$ and $e'$ of $C(U)$ are locally isomorphic. This means
that there exists an open covering $(U_i)_{i\in I}$ of $U$ such
that the restrictions  $e_i$ and $e'_i$ of respectively $e$ and
$e'$ to $U_i$ are isomorphic.

\bigskip

We say that a $2$-gerbed is bounded by the sheaf of abelian groups
$L$, if the following two conditions are satisfied:

4. Any pair  of $1$-arrows can be joined by a $2$-arrow.

 5. Let $e_U$ and $e'_U$ be a pair of objects of $C(U)$.
 For any  $1$-arrow $h:e_U\rightarrow e'_U$, there is a specified isomorphism
$L(U)\rightarrow Aut(h)$, compatible with compositions and with
restrictions $\bullet$. We say that the sheaf $L$ is  the band of
the $2$-gerbe, or that the gerbe is bounded by $L$.

\bigskip

{\bf 7.2. Classifying cocycle of a $2$-gerbe.}

\bigskip

Let $(U_i)_{i\in I}$ be an open covering  of $N$ such that
$C(U_i)$ is not empty. Consider an object $e_i$ of $C(U_i)$, and
$g_{ij}:r_{U_{ij},U_j}(e_j)\longrightarrow r_{U_{ij},U_i}(e_i)$,
there exists a $2$-arrow $h_{i_1i_2i_3}:{g_{i_1i_2}}^{i_3}
{g_{i_2i_3}}^{i_1}\rightarrow {g_{i_1i_3}}^{i_2}$, and on
$U_{i_1i_2i_3i_4}$ a $2$-arrow $u_{i_1i_2i_3i_4}$ which verifies:

$$
{h_{i_1i_2i_4}}^{i_3}(Id\circ
{h_{i_2i_3i_4}}^{i_1})=u_{i_1i_2i_3i_4}({h_{i_1i_3i_4}}^{i_2}(
{h_{i_1i_2i_3}}^{i_4}\circ Id)).
$$

The family $u_{i_1i_2i_3i_4}$ is the classifying $2$-cocycle of
$C$, if the sheaf $L$ is commutative, it is a Cech cocycle in the
classical sense, and the set of isomorphic classes of $2$-gerbes
bounded by $L$ is isomorphic to $H^3(N,L)$. If $L$ is not
commutative, we define $H^3(N,L)$ to be the set of isomorphic
classes of $2$-gerbes bounded by $L$.

\bigskip

In [20] we have given a simplified version of $2$-gerbes, that we
have named $2$-gerbed towers.

\bigskip

{\bf Definition 7.2.1.}

A {\bf $2$-gerbed tower} defined on $N$, is  defined by a gerbe
$C$ on $N$ and for each object $e_U$ of $C(U)$, a gerbe $C_1(e_U)$
defined on $U$ such that the following conditions are satisfied:

\bigskip

 (i) For each embedding map $U\rightarrow V$, there
 exists a restriction functor $r^1_{U,V}:C_1(e_U)\rightarrow
 C_1(r_{U,V}(e_U))$ such that $r^1_{V,W}\circ r^1_{U,V}=r^1_{U,W}$,
 where $r$ is the restriction functor of the gerbe $C$.

\bigskip

 (iii) There exits a commutative sheaf $L_1$ defined on $N$, such that for each
 object $e_U$ of $C(U)$, the band of  $C_1(e_U)$ is the
 restriction of $L_1$ to $U$.

\bigskip

(ii) For each morphism $h:e_U\rightarrow e'_U$ of objects of
$C(U)$, there exists a functor $h^*:C_1(e_U)\rightarrow C_1(e'_U)$
which is compatible with restrictions, such that for a morphism
$h':e'_U\rightarrow e"_U$, there exists a natural transformations
between the functor $(h'h)^*$ and ${h'}^*{h}^*$. We suppose also
the functors $(h'h)^*$ and ${h'}^*{h}^*$ coincide on objects. This
implies the existence of an element $l_{h',h}$ of $L_1(U)$ such
that $(h'h)^*=l_{h',h}\circ {h'}^*h^*$.

\bigskip

{\bf 7.3. The classifying cocycle of a $2$-gerbed tower.}

\bigskip

We can associate to a $2$-gerbed tower, a $3$-Cech cocycle defined
as follows: Consider an object $e_i$ of $C(U_i)$ and a morphism
$g_{ij}:r_{U_{ij},U_j}(e_j)\rightarrow r_{U_{ij},U_i}(e_i)$. The
arrow
$c_{i_1i_2i_3}={g_{i_3i_1}}^{i_2}{g_{i_1i_2}}^{i_3}{g_{i_2i_3}}^{i_1}$
is  the Cech classifying cocycle of the gerbe  $C$. It can be
identified to an element of the band of $C$.

 The classifying  cocycle of the $2$-gerbed tower is defined by considering  the family of
automorphisms
$$
c_{i_1i_2i_3i_4}=({c_{i_2i_3i_4}}^{i_1})^*({-c_{i_1i_3i_4}}^{i_2})^*
({c_{i_1i_2i_4}}^{i_3})^*({-c_{i_1i_2i_3}}^{i_4})^*
$$

\bigskip

Property $(iii)$ implies that $c_{i_1i_2i_3i_4}$ is an element of
$L_1(U_{i_1}..{i_4})$. Contrary to the case of $2$-gerbes, it is
 after having defined the classifying cocycle, that we set the
axiom concerning the gluing property for objects:

\bigskip

{\bf Gluing property for objects.}

\bigskip

Suppose that the cohomology class of the classifying cocycle of a
$2$-gerbed tower is zero. Let  $(U_i)_{i\in I}$ be the open
covering of $N$ used to construct the cocycle. Then there exists a
gerbe $C_0$, such that for each open subset $U$ of $N$, the
restriction of $C_0$ to $U\cap U_i$ is $C_1({e_i}_U)$, where
${e_i}_U$ is the restriction of $e_i$ to $U_i\cap U$, and $e_i$ is
the object of $C(U_i)$ used to construct the $2$-cocycle.

\bigskip

{\bf Proposition 7.3.1.}

{\it Let $(C,C_1)$ be a $2$-gerbed tower defined on $N$, the
correspondence defined on the category of open subsets of $N$ as
follows:

 To each open set $U$ of $N$, $C'(U)$ is the bicategory whose
objects are gerbes $E_U$ such that  for every open covering
$(U_i)_{i\in I}$ of $U$ such that $C(U_i)$ is not empty,
 the restriction of $E_U$ to $U_i$
 is isomorphic to a gerbe $C_1(e_i)$ where $e_i$ is an object of $C(U_i)$.
  A $1$-arrow $h:E_U=C_1(e_U)\rightarrow E'_U=C_1(e'_U)$ between
two objects of $C'(U)$ is a functor $h^*$, where $h:e_U\rightarrow
e'_U$ is an arrow. A $2$-arrow is a natural transformation $l_U$
between two $1$-arrows which coincide on objects.}

\bigskip

{\bf Proof.}

\bigskip

Gluing property for objects.

\bigskip

Consider an open covering family $(U_i)_{i\in I}$ of $N$. Let
$E_i=C_1(e_i)$ be an object of $C'(U_i)$, a morphism between
$g_{ij}:E^i_j\rightarrow E^j_i$ is a functor
${h_{ij}}^*:C_1(e^i_j)\rightarrow C_1(e^j_i)$, where
$h_{ij}:e^i_j\rightarrow e^i_j$ is an arrow. A $2$-arrow between
${{h_{i_1i_2}}^{i_3}}^*{{h_{i_2i_3}}^{i_1}}^*$ and
${{h_{i_1i_3}}^{i_2}}^*$, is a natural transformation

$$
c_{i_1i_2i_3}:{{h_{i_1i_2}}^{i_3}}^*{{h_{i_2i_3}}^{i_1}}^*\rightarrow
{{h_{i_1i_3}}^{i_2}}^*,
$$
defined by  an element of $L_1(U_{i_1i_2i_3})$. The fact that:

$$
{c_{i_1i_3i_4}}^{i_2}({c_{i_1i_2i_3}}^{i_4}\circ
Id)={c_{i_1i_2i_4}}^{i_3}(Id\circ {c_{i_2i_3i_4}}^{i_1})
$$

is equivalent to the gluing property of objects of a $2$-gerbed
tower. This implies by definition the existence of an object $E_U$
whose restriction to $U_i$ is $E_i$.

\bigskip

{ Gluing conditions of arrows.}

\bigskip

Let $E_U$ and $E'_U$ be two respective objects of $C'(U)$. The
correspondence defined on the category of open subsets of $U$ by
$V\rightarrow Hom({E_U}_{\mid V},{E'_U}_{\mid V})$, where
${E_U}_{\mid V}$ and ${E'_U}_{\mid V}$ are the respective
restrictions of $E_U$ and $E'_U$ to $V$ is a sheaf of categories
since it is the sheaf of morphisms between two gerbes.

\bigskip

Let $U$ be an open subset of $N$, the objects $E_U$, and $E'_U$ of
$C'(U)$ are locally isomorphic, since the restrictions of $E_U$
and $E'_U$ to an open cover of $(U_i)_{i\in I}$ of $U$ such that
the objects
 of $C(U_i)$ are isomorphic  are isomorphic. If we replace $U$ by $N$, and choose a covering family
such that $C(U_i)$ is not empty, we obtain that $C'(U_i)$ is not
empty.

The set of automorphisms of a $1$-arrow  is isomorphic to $L_1(U)$
by definition $\bullet$

\bigskip

The notion of $2$-gerbed tower is easier to understand than the
one of $2$-gerbe, principally because, we do not need the notion
of bicategory to define it. In practice, many the examples of
$2$-gerbed are defined using the notion of $2$-gerbed tower,
another advantage of this notion is the fact that the classifying
cocycle of a $2$-gerbed tower $(C,C_1)$ is the image of a
$2$-cocycle, that is, the classifying cocycle of $C$ by a
connecting  morphism in cohomology.

\bigskip

{\bf 8. The general case.}

\bigskip

We will now describe the $2$-gerbe and $2$-gerbed towers bounded
by the sheaf of locally constant ${\R}$-functions which represent
the geometric obstruction to extend $\omega$ to $P$ when the
cohomology class $[\omega]$ of $\omega$ is not necessarily
integral.

\medskip

 Let $U$ be an open subset of $N$, and $[\Omega_U]$ an extension
 of   $[\omega]$ to $p^{-1}(U)$. We cannot define a  $T^1$-bundle
 over $p^{-1}(U)$  (as in the integral case) whose Chern class is $[\Omega_U]$.

\bigskip

{\bf Definitions 8.1.}

 We  denote by $C'_F(\Omega,p^{-1}(U))$
the gerbe defined on $p^{-1}(U)$ which is the obstruction of the
class $[\Omega_U]$ to be trivial. see 2.4 $\bullet$

\bigskip

Let $U$ be an open subset of $N$, we define the bicategory
$C^2_F(\omega)(p^{-1}(U))$ to be the class whose elements are
categories $C'_F(\Omega,p^{-1}(U))$. Let $e_1$ and $e_2$ be two
objects of $C^2_F(p^{-1}(U))$, a $1$-arrow
$f:e^1=C'_F(\Omega_1,p^{-1}(U))\rightarrow
e^2=C'_F(\Omega_2,p^{-1}(U))$ is an isomorphism of gerbes between
$e^1$ and $e^2$, and a $2$-arrow is a natural transformation
between those functors.

More precisely, on a contractible cover $(U'_i)_{i\in I}$ of
$p^{-1}(U)$, the restrictions of the objects of $e^1$ are torsors
whose objects are isomorphic to trivial ${\R}$-bundles $U'_i\times
{\R}$,  a $1$-arrow $f$ is defined by the respective objects
$e^1_i$ and $e^2_i$ of the respective restrictions of $e^1$ to
$U'_i$, and of $e^2$ to $U'_i$, and an isomorphism $f_i$ between
$e^1_i$ and $e^2_i$. Due to the natural properties of $f$, we can
use  these morphisms to rebuild completely $f$. This implies that
these datas satisfy the following properties:

The identification of $e^1_i$ and $e^2_i$ to $U'_i\times {\R}$
allows to represents $f_i$ has a morphism of the trivial torsor
$U'_i\times {\R}$, the fact that $f$ behave naturally in respect
with restrictions implies the existence of a morphism  $u_{ij}$ of
${U'}_{ij}\times{\R}$ such that $f_i=u_{ij}f_j$. We have
${u_{i_1i_2}}^{i_3}{u_{i_2i_3}}^{i_1}={u_{i_1i_3}}^{i_2}$. The map
$u_{ij}$ is a translation by an element of  ${\R}$. The family
$(u_{ij})_{i,j\in I}$ defines a $1$-cocycle, thus a closed
$1$-form on $p^{-1}(U)$. Conversely, a $1$-cocycle of the sheaf of
locally constant ${\R}$-maps defines a torsor, and a $1$-arrow
between $e^1$ and $e^2$ by using the previous identification of
$e^1_i$ and $e^2_i$ to $U'_i\times {\R}$.

 Using the identification above, a $2$-arrow is defined locally by a chain of constant
sections $u_i$ defined on $U'_i\times {\R}$ such that
$u_i=u_{ij}u_j$.  Thus a morphism between two objects is defined
by a $1$-cocycle of the sheaf of locally constant
${\R}-$functions, that is a torsor, and a $2$-arrow is an element
of ${\R}$.

\bigskip

{\bf Theorem 8.2.}

{\it The correspondence $p^{-1}(U)\rightarrow
C^2_F(\omega)(p^{-1}(U))$ defines a $2$-gerbe such that the
cohomology class of its  classifying cocycle, is the obstruction
for extending $[\omega]$ to $P$.}

\bigskip

{\bf Proof.}

\bigskip

Gluing conditions for objects.

\bigskip

Consider an open covering  $(U_i)_{i\in I}$ of an open subset $U$
of $N$, and $(e_i,[\Omega_i])$ an object of
$C^2_F(\omega)(p^{-1}(U_i))$ where $[\Omega_i]$ is a cohomology
class defined on $p^{-1}(U_i)$ which extends $[\omega]$. Suppose
that there exists $1$-arrows $h_{ij}:e^i_j\rightarrow e^j_i$, a
$2$-arrow
$d_{i_1i_2i_3}:{h_{i_1i_2}}^{i_3}{h_{i_2i_3}}^{i_1}\rightarrow
{h_{i_1i_3}}^{i_2}$ such that
${d_{i_1i_3i_4}}^{i_2}({d_{i_1i_2i_3}}^{i_4}\circ
Id)={d_{i_1i_2i_4}}^{i_3}(Id\circ {d_{i_2i_3i_4}}^{i_1})$.  The
maps $d_{i_1i_2i_3}$ can be identified with a $2$-Cech cocycle of
the sheaf of locally constant  ${\R}$-functions defined on
$p^{-1}(U)$. We can identify it using the De Rham Weil isomorphism
with an element $[\Omega_U]$ of $H^2(p^{-1}(U),{\R})$. The class
$[\Omega_U]$ is the classifying cocycle of a gerbe $e_U$ defined
on $U$ bounded by the sheaf of locally constant ${\R}$-functions.
We have to show now that this gerbe is an element of
$C^2_F(p^{-1}(U))$.

The fact that the family $d_{i_1i_2i_3}$ is a $2$-Cech cocycle,
implies that there exists a gerbe bounded by the sheaf of locally
constant functions whose restriction to $U_i$ is $e_i$, (See the
proof of the classifying theorem for gerbes presented in the book
of Breen [2]) this gerbe is isomorphic to $e_U$. This implies that
 the restriction of $[\Omega_U]$ to
$p^{-1}(U_i)$ is the classifying cocycle of $e_i$, and that
$[\Omega_U]$ extends $[\omega]$, since the restriction of
$[\Omega_U]$  to $U_i$ is $[\Omega_i]$. We deduce that it is the
cohomology class of the classifying cocycle of an element of
$C^2_F(\omega)(p^{-1}(U))$ whose restriction to $p^{-1}(U_i)$ is
$e_i$.

\bigskip

Gluing conditions for arrows.

\bigskip

Let $e$, and $e'$ be a pair of objects of
$C^2_F(\omega)(p^{-1}(U))$, the correspondence defined on the
category of open subsets of $U$ by $U'\rightarrow Hom(e_{\mid
U'},e'_{\mid U'})$ is a sheaf of categories, since  it is the
sheaf of categories of morphisms between two gerbes.

This shows that $C^2_F(\omega)$ is a sheaf of $2$-categories.

\bigskip

 Consider
an open covering $(U_i)_{i\in I}$ of $N$ by contractible open
sets, since $H^*(U_i\times F)=H^*(F)$, we can extend $[\omega]$ to
$p^{-1}(U_i)=U_i\times F$, and two such extension classes are
equal to the class $[\omega]$ as shows the  identification
$H^*(U_i\times F)=H^*(F)$. This implies that
$C^2_F(\omega)(p^{-1}(U_i))$ is not empty, an its objects are
isomorphic.

The sheaf of $2$-categories $C^2_F(\omega)$ is bounded by the
sheaf of ${\R}$-locally constant functions defined on $P$. This is
shown in the paragraph above this theorem.

Suppose that the class of the classifying cocycle of
$C^2_F(\omega)$ vanishes, then the $2$-gerbe has a global section
$e$, its restriction to $p^{-1}(U_i)$ is an element of
$C^2_F(\omega)(p^{-1}(U_i))$ whose classifing cocycle extends
$[\omega]$. This implies that the classifying cocycle of $e$
extends $[\omega]$ $\bullet$

\bigskip

The cocycle defined by McDuff, and the classifying cocycle $c^2_F$
of the $2$-gerbe $C^2_F(\omega)$ solve the same geometric problem:
decide if the class $[\omega]$ can be extended to the total space
of the symplectic bundle $p:P\rightarrow N$ whose typical fiber is
$(F,\omega)$. We will show now that they are related by an
isomorphism of cohomology groups.

 Suppose that the family $(U_i)_{i\in I},g_{ij}$ defines the
coordinate changes of $P$. Let $\hat Symp(F,\omega)$ be the
universal cover of $Symp(F,\omega)$. Consider an element $h_{ij}$
of $Ham^s$ such that $g_{ij}(x)h_{ij}=g'_{ij}(x)$ is contained in
$Symp(F,\omega)_0$, and a lift: $\hat{g'_{ij}}:U_i\cap
U_j\rightarrow \hat Symp(F,\omega)$ of the functions $g'_{ij}$.
Remark that an element $\hat {g'_{ij}}(x)$ is an equivalence class
of a path in $c:[0,1]\rightarrow Symp(F,\omega)$. We choose a path
$u_{ij}$ which represents it and set:

$$
\int_0^1 \omega({d\over {dt}}\hat {u_{ij}}(x),.)=g"_{ij}(x)
$$

\bigskip

{\bf Proposition. 8.3}

{\it The chain
$c_{i_1i_2i_3}={g"_{i_2i_3}}^{i_1}-{g"_{i_1i_3}}^{i_2}+{g"_{i_1i_2}}^{i_3}$
is a $2$-Cech cocycle whose cohomology class is identified using
the Cech-Weil isomorphism to the McDuff obstruction class.}

\bigskip

{\bf Proof.}

The element $g"_{ij}(x)$ is a lift of $F_s(g"_{ij})$ in
$H^1(F,{\R})$, since the restriction of $F_s$ to
$Symp(F,\omega)_0$ factors by the flux. This implies that
$g"_{ij}$ represents the classifying cocycle of the
$H^1(F,{\R})/H^1(F,P_{\omega})$-bundle (see 2.6.5) whose
coordinate changes are the functions $F_s(g'_{ij})$ $\bullet$

\bigskip

We can use the Cech-Weil isomorphism to identify $c_{i_1i_2i_3}$
to a  closed $2$-form $\Omega'$ defined on $N$ which take values
in the vector bundle $p_{\omega}$ of closed $P_{\omega}$ $1$-forms
defined on $F$ induced by $g_{ij}$. Let $\Omega(F,P_{\omega})$ be
the vector space of closed $P_{\omega}$ $1$-forms defined on $F$.
The bundle $p_{\omega}$ is the quotient of the union of $U_i\times
\Omega(F,P_{\omega})$ by the following transitions functions:

$$
(x,\alpha)\rightarrow (x,{g_{ij}(x)}^*(\alpha))
$$

where ${g_{ij}(x)}^*(\alpha)(y)$ is defined by:

$$
{g_{ij}(x)}^*(\alpha)(y)=\alpha(d({g_{ij}(x)}^{-1})(y))
$$

 The identification of $c_{i_1i_2i_3}$ to $\Omega'$
 defines a $3$-form $\Omega$ on $P$ by
$\Omega(x,y,z)=\Omega'(x,y)(z)$ where $x,y$ are elements of $T_nN$
the tangent space of $N$ at $n$, and $z$ is an element of the
tangent space to the fiber at $n$.

Consider the Leray-Serre spectral sequence related to the
fibration $p:P\rightarrow N$ the McDuff obstruction class is an
element of $E_2^{2,1}$ which converges to $[c_F^2(\omega)]$.

\bigskip

{\bf Theorem 8.3.}

{\it  Under the notation just above, the cohomology class of
$\Omega$
  is the obstruction to lift
$[\omega]$ to $P$. Its cohomology class  can identified to the
class of the classifying cocycle of $C^2_F(\omega)$.}

\bigskip

{\bf Proof.}

Let $e_i$ be the gerbe defined on $U_i\times F$ whose classifying
cohomology class is the image of  the class of the $2$-form
$\Omega_i$ which is the product of $0$ and $\omega$ by the
Cech-Weil isomorphism. The gerbe $e_i$ is an object of
$C_F^2(\omega)(U_i)$. The morphism $g"_{ij}$ defined at the
proposition above is a morphism between $e_j^i$ and $e_i^j$. This
implies that $c_{i_1i_2i_3}$ represents also the classifying
cocycle of $C_F^2(\omega)$ $\bullet$

\bigskip

{\bf 8.2 Hamiltonian reduction and closed connection forms.}

\bigskip

 We have
 given a gerbe formulation to the  problem of the existence of a Hamiltonian reduction, by defining the gerbe
 $C^1_F(\omega)$, now we are going to show how the classifying cocycle of
 $C^1_F(\omega)$ and $C^2_F(\omega)$ are related.

\bigskip

The link between the classifying cocycles of $C^1_F(\omega)$ and
$C^2_F(\omega)$ appears clearly by considering the $2$-gerbed
towers defined as follows:

\bigskip

{\bf Definition 8.2.1.}

Consider $U$, an open set of $N$, $e_U$ an object of
$C^1_F(\omega)(U)$, it is a Hamiltonian structure defined on the
restriction of the symplectic fibration $p:P\rightarrow N$ to $U$.
We deduce that there exists an extension $[\Omega_U]$ of
$[\omega]$ to $p^{-1}(U)$ whose holonomy defines the Hamiltonian
reduction of $e_U$. Denote by $C_2(e_U)$  the gerbe which
represents the obstruction of $[\Omega_U]$ to be trivial. We have
just defined a $2$-gerbed tower $(C_F^1(\omega),C_2)$ $\bullet$

\bigskip

  Let $L$ be the band of the gerbe $C^1_F(\omega)$, and $L_0$ the sheaf
  of locally constant ${\R}$-functions defined on $P$. We define the following sheaf
   $L'$ on $P$:  suppose that $e_U$ is an object of
$C_F^1(\omega)(U)$,  $V$ an open subset of $p^{-1}(U)$, and $e_V$
 an object of $C_2(e_U)(V)$. An automorphism $g$ of $e_U$ map
$e_V$ to the object ${g^{-1}}^*(e_V)$ of $C_2(e_U)(g(V))$, given $
c\in {\R}$, for each morphism $h:e_V\rightarrow e'_V$ between
objects of $C_2(e_U)(V)$ we consider the morphism between
${g^{-1}}^*(e_V)\rightarrow {g^{-1}}^*(e'_V)$ defined by composing
${g^{-1}}^*(h)$ by the translation by $c$ fiber by fiber. The
sheaf generated by the set of actions on the gerbe $C_2(e_U)$ that
we have just defined is $L'$. (It does not depend of $e_U$). We
can suppose that $L$ is defined on $P$ by setting $L(U)=L(p(U)),
U\subset P$. We have the exact sequence:

$$
1\longrightarrow L_0\longrightarrow L'\longrightarrow
L\longrightarrow 1
$$

This gives rise to the following exact sequence in cohomology:

$$
H^2(P,L_0)\longrightarrow H^2(P,L')\longrightarrow
H^2(P,L)\longrightarrow H^3(P,L_0)
$$

Here if $E$ is a sheaf defined on $P$, the space $H^2(P,E)$ is the
space of isomorphism classes of gerbes bounded by $E$ defined on
$P$. The space $H^3(P,E)$ is the space   of isomorphism classes of
$2$-gerbes bounded by $E$. See Breen [2].

The next result show that the class of the classifying cocycle of
the $2$-gerbe tower $(C^1_F(\omega),C_2)$ is the image of the
class of the classifying cocycle of $C^1_F(\omega)$ by the map
$H^2(P,L)\longrightarrow H^3(P,L_0)$.

\bigskip

{\bf Proposition 8.2.2.}

{\it The class of the classifying cocycle $c_F^2(\omega)$ of
$C^2_F(\omega)$, is the image of the class of the classifying
cocycle $c^1_F(\omega)$ of $C^1_F(\omega)$, by the map
$H^2(P,L)\rightarrow H^3(P,L_0)$. Suppose that there exists a
Hamiltonian reduction of the bundle $P\rightarrow N$, then we can
extend $[\omega]$ to $P$.}

\bigskip

{\bf Proof.}

The classifying cocycle of this $2$-gerbed tower is defined as
follows, consider an object $e_i$ of $C^1_F(\omega)(U_i)$, and a
map $u_{ij}:e^i_j\rightarrow e^j_i$, the map
$c_{i_1i_2i_3}={u_{i_3i_1}}^{i_2}{u_{i_1i_2}}^{i_3}{u_{i_2i_3}}^{i_1}$
is an automorphism of ${e^{i_1i_2}}_{i_3}$, we can lift it to a
map ${c^*}_{i_1i_2i_3}$ of $C_2({e^{i_1i_2}}_{i_3})$, the Cech
boundary $c_{i_1i_2i_3i_4}$ of the chain ${c^*}_{i_1i_2i_3}$ is
the classifying cocycle of the $2$-gerbed tower. It appears that
$c_{i_1i_2i_3i_4}$ is the image of $c_{i_1i_2i_3}$ by the
connecting  map $H^2(P,L)\rightarrow H^3(P,L)$. Considered as a
$2$-gerbe the $2$-gerbed tower involved here is a subgerbe of
$C^2_F(\omega)$, since for each object $e_U$ of
$C_F^1(\omega)(U)$, the gerbe $C_2(e_U)$ is an object of
$C_F^2(\omega)(U)$.  This shows that if $[c^1_F(\omega)]$
vanishes, then $[c^2_F(\omega)]$ also vanishes $\bullet$

\bigskip

This result is shown by McDuff in [14] by using the
Guillemin-Lerman-Sternberg construction.

\bigskip

\bigskip

{\bf 9. Quantization of the symplectic gerbe.}

\bigskip

Let $(F,\omega)$ be a symplectic manifold, when the class
$[\omega]$ is  integral, there exists a line bundle $L$ over $F$
whose chern class is $[\omega]$. This line bundle is endowed with
a hermitian metric. The hermitian space of sections
$L^2(F)=\{u:F\rightarrow L: \int_F\mid u\mid^2<+\infty\}$ is the
quantization of the manifold. The elements of this Banach space
are used in theoretical physic, to describe evolution of
particles.

The goal of this part is to associate to any symplectic form, a
hermitian space endowed with a Hermitian form, which is a
candidate to represents the phase space  in quantum mechanic.

\bigskip

Let $C(\omega)$ be the symplectic gerbe defined on $F$, which
represents the obstruction of $[\omega]$ to be integral see 2.4.
Consider an open covering $(U_i)_{i\in I}$ of $U$, and $e_i$ an
object of $C(\omega)(U_i)$. We can define the gerbe $L(\omega)$ on
$F$, such that $L(\omega)(U)$ is the category, whose objects are
$(e_U,{e'}_U)$ where $e_U$ is an object of $C(\omega)(U)$, and
$e'$ the ${\C}$-line vector bundle over $U$, whose transition
functions are the transition functions of $e_U$. The object of
$L(\omega)(U)$ are endowed with a canonical connective structure
$Co$ see 2.4.  An element of $Co((e_U,e'_U))$ is a connection on
$e_U$ whose curvature is the restriction of $\omega$ to $U$.
 A morphism between two objects $(e_U,{e'}_U)$ and ${e^1}_U),{e'}^1_U)$
 of $L(\omega)(U)$ is a morphism $e_U\rightarrow {e^1}_U$.
 The correspondence defined on the category of open subsets of
$F$ by $U\rightarrow L(\omega)(U)$ is a gerbe.

\bigskip

To perform the quantization we need to define the notion of
sections. We will propose this definition of sections of vectorial
gerbes.

\bigskip

{\bf Definition 9.1.}

Let $(U_i)_{i\in I}$ be an open  covering family of $F$, such that
$L(\omega)(U_i)$ is not empty, and $(e_i,e'_i)$  an object
$L(\omega)(U_i)$. A section $u$ of $(e'_i)_{i\in I}$ is a family
of sections $u_i:U_i\rightarrow e'_i$ such that the union of
supports of $u_i$ is compact $\bullet$

\bigskip

 We denote by $V((e_i)_{i\in I})$ the vector space generated by
those sections of $(e'_i)_{i\in I}$. This   vector space  is
endowed with a Hermitian metric defined by

$$
<u,v>=\sum_{i\in I}\int_{{e'_i}}<u_i,v_i>_{e'_i}
$$

 where $<,>_{{e'}_i}$ is the Hermitian metric of $e'_i$.

\bigskip

 For each function $f$, and each section $(u_i)_{i\in I}$. We can
define

$$
L_f(u_i)={\nabla_{e'_i}}_{X_f}u_i +2i\pi fu_i,
$$

 where $X_f$ is the
Hamiltonian of $f$, and $\nabla_{e'_i}$ a connection defined on
$e'_i$ whose curvature is the restriction of $\omega$ to $U_i$.
The vector field $X_f$ is the vector field such that
$\omega(X_f,.)=-df$.

\bigskip

{\bf Proposition 9.2.}

{\it The family of $L_f(u_i)$ defined a section $L_f(u)$. The map

$$
f\rightarrow L_f
$$

verifies

$$
[L_f,L_g]=L_{\{f,g\}}
$$}

\bigskip

{\bf Proof.}

We have to show that $L_f(u_i)$ has a compact support, and that
the union of support of the family $(L_f(u_i))_{i\in I}$ is
compact. The support of $fu_i$ and ${\nabla_{e'_i}}_{X_f}(u_i)$
are contained in the support of $u_i$.
 The fact that $
[L_f,L_g]=L_{\{f,g\}}$ is classical. $\bullet$

We have obtained a Souriau-Kostant quantization.

\bigskip

We can define using the classifying theorem of Giraud [6] the
gerbe $L'(\omega)$ on $F$, such that $L'(\omega)(U)$ is a set of
flat ${\C}$-bundles defined on $U$, and the cohomology class of
the classifying cocycle of $L'(\omega)$ is the obstruction of the
class $[\omega]$ to be integral. This construction of this gerbe
using [3] shows that this gerbe is flat, the objects of
$L'(\omega)(U)$ are locally flat ${\C}$-bundles, and morphisms are
morphisms of locally flat ${\C}$-bundles.

\bigskip

For $L'(\omega)$, we can also define the following space of
sections. Consider an open covering $(U_i)_{i\in I}$ of $F$ by
contractible subsets, $e'_i$ an object of $L'(\omega)(U_i)$,
$g_{ij}:{e'}_j^i\rightarrow {e'}^j_i$ a family of isomorphisms. A
section $u=(u_i)_{i\in I}$ is a family of sections
$u_i:U_i\rightarrow e'_i$ such that $u_i=g_{ij}(u_j)$. We denote
by $V(e_i,g_{ij})$ the set of those sections. It is a vector space
which can be endowed with the following scalar product.

 Consider a
partition of the unity $p_i$ subordinate to $(U_i)_{i\in I}$. Let
$u=(u_i)_{i\in I}$, and $u'=({u}'_i)_{i\in I}$ be sections of
$V(e_i,u_{ij}))$. We set

$$
<u,v>\sum_{i\in I} \int<p_iu_i,p_i{u'}_i>
$$

For each differentiable function $f$ defined on $F$ we can define
the operator $L_f$ which acts on the section $u=(u_i)_{i\in I}$
by:

$$
L_f(u_i)=\nabla_{X_f}u_i+ 2i\pi fu_i
$$

The operator $L_f$ is well defined. Since on $U_i\cap U_j$, we
have $L_f(u_i)=u_{ij}L_f(u_j)$ since the gerbe $C(\omega)$ is
flat, and the map $u_{ij}$ are identified using a trivialization
with the multiplication by an element of $T^1$ in the trivial
bundle $U_i\cap U_j\times {\C}$.

\bigskip

{\bf  Quantization of other structures.}

\bigskip

The methods of quantization of Kostant-Souriau have been extended
in many directions. Here we present a quantization described in
[19]

\bigskip

Consider a manifold $M$, such that the ring $C^{\infty}(M)$ of
differentiable functions of $M$ is endowed with a bracket:

$$
\{,\}:C^{\infty}(M)\times C^{\infty}(M)\longrightarrow
C^{\infty}(M)
$$

such that $(C^{\infty}(M),\{,\})$ is a Lie algebra  and there is a
${\R}$-linear map:

$$
H:C^{\infty}(M)\longrightarrow \chi(M)
$$

$$
f\rightarrow X_f
$$

where $\chi(M)$ is the space of vector fields of $M$, such that

$$
X_{\{f,g\}}=[X_f,X_g]
$$

The map

$$
C^{\infty}(M)\longrightarrow End(C^{\infty}(M))
$$

$$
f\longrightarrow (g\rightarrow X_f(g))
$$

is a representation of the Lie algebra $C^{\infty}(M)$. We denote
$H_C^*(M)$ the cohomology of this representation. The
correspondence:

$$
C^{\infty}(M)\times C^{\infty}(M)\longrightarrow C^{\infty}(M)
$$

$$
\Lambda_M(f,g)\longrightarrow X_f(g)-X_g(f)-\{f,g\}
$$

is a $2$-cocycle of this representation.

\bigskip

There is a canonical map $C':H^*_{De Rham}(M)\rightarrow
H^*_C(M)$. Defined on a chain by
$C'(h)(f_1,...,f_p)=h(X_{f_1},...,X_{f_p})$. In [19], it is shown
that if there exists a line bundle $L\rightarrow M$ such that
$C'(\Omega)=\Lambda_M$, then the structure is quantizable: that is
there exists a representation:

$$
P:C^{\infty}(M)\longrightarrow End(L^2(L))
$$

which verifies

$$
P(\{f,g\})=[P_f,P_g]
$$

$$
P(f)=\nabla_{X_f}+2i\pi f
$$

where $\nabla$ is the hermitian connection of the bundle.

\bigskip

Let $(U_i)_{i\in I}$ be a contractible open covering of $M$ by
charts. We can restrict the bracket $\{,\}$ to $U_i$. Suppose that
on $2$-chains, the map $C$ restricted to $U_i$ is surjective on
closed forms. that is there exists a $2$-closed form
$\Omega_{U_i}$ on $U_i$ such that $C(\Omega_{U_i})=\Lambda_{U_i}$.
The form $\Omega_{U_i}$ is the Chern class of a connection defined
on $U_i\times {\C}$.

\bigskip

We can define on $M$ the gerbe $D$, such that for each open set
$U$ of $M$, $D(U)$ is the category of line bundles over $U$
endowed with a connection whose curvature $\Omega_U$ verifies:

$$
C(\Omega_U)=\Lambda_U
$$

Let $e_i$ be an object of $D(U_i)$, we consider the family
$(u_i)_{i\in I}$, where $u_i:U_i\rightarrow e_i$ is a section of
$e_i$, whose support is compact, and the union of support of $u_i$
is compact. The family of $(u_i)_{i\in I}$ is a Hermitian space.
On $e_i$ we consider the   connection $\nabla_{e_i}$ whose
curvature is the restriction of $\Omega_{U_i}$

The representation

$$
f\rightarrow {\nabla_{e_i}}_{X_f}+ 2i\pi f
$$

defines a quantization of $M$.

\bigskip

\bigskip
\bigskip

\centerline{\bf References.}

\bigskip

[1] A. Banyaga. Sur la structure du groupe des automorphismes
symplectiques qui preservent une forme symplectique.
Communicationes Mathematicae Helveticae. 53 1978. 174-227

\medskip

 [2] L, Breen. On the classification of $2-$gerbes and $2-$stacks.
Asterisque, 225 1994.

\medskip

 [3] J.L Brylinski, Loops spaces, Characteristic Classes and
Geometric Quantization, Progr. Math. 107, Birkhauser, 1993.

\medskip

 [4]  J.L Brylinski, D.A McLaughlin, The geometry of degree four
characteristic classes and of line bundles on loop spaces I. Duke
Math. Journal. 75 (1994) 603-637.

\medskip

[5] P. Deligne. Theorie de Hodge III. Institut des Hautes Etudes
Scientifiques Publications 44 5-77.

\medskip

[6] J. Giraud. Cohomologie non abelienne.

\medskip

[7] Godement. Theories des faisceaux

\medskip

[8] M. Gotay, R. Lashof, J. Sniatycki, A. Weinstein. Closed form
on symplectic bundles. Commentarii. Mathematicae. Helvetici 58
(1983) 617-621.

\medskip

[9] V. Guillemin, E. Lerman, S. Sternberg. Symplectic fibrations
and multiplicity diagrams. Cambridge University Press 1998

\medskip

[10] A. Grothendieck. Cohomologie Etale. Seminaire de Geometrie
algebrique. 4.

\medskip

[11] J. Kedra, D McDuff. Homotopy properties of Hamiltonian group
actions. SG/0404539

\medskip

[12] J. Leslie. On a differentiable structure for he group of
diffeomorphisms.  Topology 6. 1967 263-271

\medskip

[13] F. Lalonde, D. McDuff. Symplectic structures on fiber
bundles. Topology (42) 2003 309-347

\medskip

[14]. D. McDuff. Enlarging the Hamiltonian group

\medskip

[15] D. McDuff. Quantum homology of fibrations over $S^2$.
International. Math. Journal. 11 665-721 2000

\medskip

 [16] D. McDuff, D. Salamon. Introduction to symplectic topology.
O.U.P

\medskip

[17] D. McDuff, D. Salamon. J-holomorphic curves and quantum
cohomology. Amer. Math. Soc. Lectures Notes 6 Providence 1994

\medskip

 [18] Maakay, Picken,  Holonomy an parallel transport for gerbes.
Advances in Math.

\medskip

 [19] Leon, Marrero, Padron, On the geometric prequantization of
brackets. Rev. Acad. Cien. Serie. Math. vol 95.

\medskip

 [20] A. Tsemo. Non abelian cohomology: the point of view of gerbed
tower. preprint available at arxiv.org

[21] A. Tsemo. Fibr\'es affines. Michigan. Math. Journal. 49

\end{document}